\documentclass[11pt, leqno]{amsart}
\usepackage{amsthm, amsmath}
\usepackage{amssymb} 
\usepackage{amscd}
\usepackage[curve,matrix,arrow,frame,tips]{xy}
\usepackage{verbatim}
\usepackage{graphicx}

\newtheorem{thm}{Theorem}[section]
\newtheorem{prop}[thm]{Proposition}
\newtheorem{lemma}[thm]{Lemma}

\newtheorem{Remark}[thm]{Remark}

\newtheorem{cor}[thm]{Corollary}

\newcounter{ex}[section]

\newcommand{\cal}{\mathcal}

\newcommand{\A}{{\mathcal A}}

\newcommand{\Q}{{\bf Q}}

\newcommand{\Z}{{\bf Z}}

 \renewcommand{\O}{{\mathcal O}}

\newcommand{\Gr}{{\rm G}}

\newcommand{\Kr}{{\rm K}}

\newcommand{\GL}{{\rm GL}}

\newcommand{\SK}{{\rm SK}}

\renewcommand{\mod}{\ {\rm mod}\ }

\def\thfill{\null\nobreak\hfill}

\def\endproof{\thfill\vbox{\hrule
  \hbox{\vrule\hbox to 5pt{\vbox to 5pt{\vfil}\hfil}\vrule}\hrule}}

\renewcommand{\P}{{\cal P}}

\newcommand{\Gal}{{\rm Gal}}

\begin{document}
 
\title[$\Kr_1$ of a $p$-adic group ring]{$\Kr_{1}$ of a $p$-adic group ring I. \\ The determinantal image. }

\author[T. Chinburg]{T. Chinburg}\thanks{Chinburg is supported by NSF Grant \# DMS08-01030}
\address{Ted Chinburg, Dept. of Math\\Univ. of Penn.\\Phila. PA. 19104, U.S.A.}
\email{ted@math.upenn.edu}

\author[G. Pappas]{G. Pappas}
\thanks{Pappas is supported by NSF Grant \# DMS08-02686}
\address{Georgios Pappas, Dept. of Math\\Michigan State Univ.\\East Lansing, MI 48824, U.S.A.  }
\email{pappas@math.msu.edu}

\author[M. Taylor]{M. J. Taylor}
\thanks{ Taylor was in receipt of a Royal Society Wolfson Merit award for the majority of the period when this article was being written.}
\address{Martin J. Taylor, School of Math.\\ Univ. of Manchester\\ Manchester, M60 1QD, U.K.}
\email{Martin.Taylor@manchester.ac.uk}

\date{\today}

\maketitle


 \bigskip

\noindent {\bf Abstract:}
  {\sl  We study the $K$-group  $K_{1}$ of the group ring of a finite
group over a coefficient ring which is $p$-adically complete and admits a
lift of Frobenius.  In this paper, we consider  the  image of $K_1$ under the determinant map;
the central tool  is the group logarithm which we can define using the Frobenius lift.
Using this we prove a fixed point theorem for the determinantal image
of $K_{1}$.} 
\medskip


\bigskip

\centerline{\sc Contents}
\medskip

\noindent \S  1. Introduction\\
\S 2. Generalities  \\
\S 3. The group logarithm and determinants for $p$-groups  \\
\S 4. Descent and fixed points for $p$-groups \\
\S 5. Character action and reduction to elementary groups  \\
\S 6. Determinants for elementary groups  \\
\S 7.  Appendix: Torsion determinants 
 \bigskip
\medskip

\section{Introduction}

In this article we begin a two part study of the $\Kr$-group $\Kr_{1}(R[G])$ of a group ring $R[G]$ for a finite group $G$, 
where the coefficient ring $R$ is a $p$-adically complete ring. The study of $\Kr_1$ for such group rings, in the
case where $R$ is a $p$-adic ring of integers, played a crucial role in
understanding the Galois structure of rings of integers of tamely ramified
extensions of number fields, and, in particular, in the solution of the Fr\"{o}hlich conjecture. 
More recently there has been considerable interest in 
$\Kr_{1}$ of group rings with higher dimensional rings of coefficients in
equivariant Iwasawa theory (see for instance [FK], [K] and [RW]). The
present paper, however, is principally oriented towards the study of the
equivariant Euler characteristics of torsors of finite groups over
arithmetic schemes. The use of the results of this paper in this context
will be discussed in more detail towards the end of this introductory
section.

Recall that, for an arbitrary (unital) ring $S$ by definition we know that $\Kr_1(S)$ sits in an exact sequence
\begin{equation*}
1\rightarrow {\rm E}(S) \rightarrow \GL(S) \rightarrow
\Kr_1(S) \rightarrow 1
\end{equation*}
where $\GL(S) =\displaystyle{\lim_{\rightarrow }}\, \GL_{n}(S)$ denotes
the general linear group of $S$ and ${\rm E}(S)$ denotes the subgroup
of elementary matrices over $S$.

Throughout this paper $R$ will always denote an integral domain of
characteristic zero with field of fractions $N$ and $N^{c}$ will denote a
chosen algebraic closure of $N$. We then have a map, which we denote $
\mathrm{Det}$
\begin{equation}\label{eq1}
\mathrm{Det}:\Kr_{1}(R[G]) \rightarrow \Kr_{1} (
N^{c}[G]) =\prod_{\chi} N^{c\times }
\end{equation}
where the direct product extends over the irreducible $N^{c}$-valued characters
of $G$. We write $\SK_{1}(R[G]) =\ker  ( \mathrm{Det} )$, 
so that we have the exact sequence
\begin{equation}\label{eq2}
1\rightarrow \SK_{1}(R[G]) \rightarrow \Kr_{1}(R[G]) \rightarrow \mathrm{Det}(\Kr_1(R[G])) \rightarrow 1.
\end{equation}
We shall seek to understand $\Kr_{1}(R[G]) $ by
firstly studying $\mathrm{Det} ( \Kr_{1}(R[G])) $ in this paper (see Theorems \ref{thm3}, \ref{thm5} and \ref{thm6} below); then in [CPT1] we
shall also use the results of this paper to study $\SK_{1}(R[G])$.   Our results for $\mathrm{Det}(\Kr_1(R[G]))$ will be obtained by a generalization of the group
logarithm, as developed in [T] (see additionally [CR2] and [F]). See also
R.~Oliver's sequence of papers [O1-4], which also contains his own
independent description of the group logarithm.

  We let $\mathrm{Det}(\GL(R[G]))
 $ denote the image of $\GL(R[G])$ under
the composition of surjective maps 
\begin{equation*}
\GL(R[G])\rightarrow \Kr_{1}(R[G]) \rightarrow \mathrm{Det} (\Kr_{1}(R[G]))
\end{equation*}
so that by definition $\mathrm{Det}(\GL(R[G]))=\mathrm{Det} (\Kr_{1}(R[G]))$.
 We let $\mathrm{Det}( R [ G ] ^{\times } ) $ denote
the subgroup of $\mathrm{Det} ( \GL(R [ G]))
  $ given by the image, under $\mathrm{Det}$, of the subgroup $R[G]^{\times }=\GL_{1}(R[G])$ 
of $\GL(R[G])$ where as usual we view $R[G]^{\times }$ as embedded in $\GL(R[G])) $ by
mapping such elements into the leading diagonal position.
\medskip

\noindent \textbf{Hypothesis:}   Let $p$ be a prime number. Throughout this paper,
unless explicitly indicated   to the contrary, we shall assume that: 

(i) $R$ is
an integral domain which is torsion free as an abelian group,

(ii)  the natural map $R\rightarrow \displaystyle{\lim_{\leftarrow }}\, R/p^{n}R$ is an isomorphism, so that $R$ is 
$p$-adically complete,

(iii)   $R$ supports a lift of Frobenius, that is to
say an endomorphism $F=F_{R}:R\rightarrow R$ with
the property that for all $r\in R$ 
\begin{equation*}
F(r) \equiv r^{p}\, \mod\, pR.
\end{equation*}

Unless stated to the contrary, we shall extend $F$ to an $R$-module
endomorphism of $R[G] $ by setting $F(g) =g$ for all 
$g\in G$.

\begin{Remark}{\rm 
a) As $R$ is $p$-adically complete, then, since $R[G]$ is a free $R$-module of finite rank, $R[G]$ is of course also $p$-adically
complete. Note that for all $s\in R$ the element $ 1-ps $ has
inverse $\sum_{n=0}^{\infty } ( ps ) ^{n}$ so that the Jacobson
radical of $R$ necessarily contains $pR$. If $\{u_{n}\}$ is a $p$-adically
convergent sequence of units converging to $r$ in $R [ G ] $, then $
r$ will be congruent to   $u_{n}\mod p$ for large $n$ and so $r$
is a unit; hence $R[G] ^{\times }$ is a $p$-adically  closed
subset of $R[G]$.

b)  If $R$ is a formally smooth Tate
algebra or dagger algebra  
then $R$ possesses a lift of Frobenius (e.g see Theorem A-1 on page 158 of [C]).
The same is true for the Witt ring 
$R=W(S)$ of a $p$-adic ring $S$. In fact, if $S$ is $p$-adically complete and separated 
the same is true for $R=W(S)$. If in fact, $S$ is a perfect integral domain
of characteristic $p$, then $R=W(S)$ satisfies our hypothesis.} 
\end{Remark}

Our first result is well known for semi-local rings (see Theorem 40.31 in
[CR2]); we obtain this result for the very different class of rings $R$
satisfying the above hypothesis:

\begin{thm}\label{thm3}
Let $R$ be as in the Hypothesis. The inclusion $R[ G ]\subset \GL ( R[ G]) $
induces an equality $$\mathrm{Det}( R[ G] ^{\times }) =
\mathrm{Det}( \GL( R[ G] )):=\mathrm{Det} (\Kr_1(R[G]))
$$ in the
following two circumstances:

(a) when $ G $ is a $p$-group;

(b) when $G$ is an arbitrary finite group and when $R$ is additionally
Noetherian and normal.
\end{thm}

In the paper [CPT1], and also in future applications, we shall be
particularly interested in the completed K-groups
\begin{equation*}
\widehat{\Kr}_{i}(R[G]) =\lim_{\leftarrow
}\,\Kr_{i}(R/p^nR[G]).
\end{equation*}%
From Proposition 1.5.1 in [FK] we know that, if all the quotient rings $R/p^{n}R$ are finite, then the natural map $\Kr_1(R[G])\rightarrow 
\widehat{\Kr}_{1}(R[G])$
is an isomorphism. (See also the work of Wall in [W1] and [W2].)
Let $\widehat{\SK}_{1}(R[G])$ denote the closure of 
$\SK_{1}(R[G]) $ in $  \widehat{\Kr}_{1}(R[G])$ and define
\begin{equation*}
\widehat{\mathrm{Det}} (\Kr_{1}(R[G]))
=\lim_{\leftarrow }\frac{\mathrm{Det} ( \GL(R[G]))}{\mathrm{Det} ( 1+p^{n}{\rm M} ( R [ G ] )) }.
\end{equation*}

\begin{prop}\label{prop4}
The natural map  $$\mathrm{Det}(\Kr_1(R[G])) \rightarrow \widehat{\mathrm{Det}}(\Kr_1(R[G]))$$ is an isomorphism and the homomorphism 
$\mathrm{Det}$ gives an exact sequence
\begin{equation}\label{eq3}
1\rightarrow \widehat{\SK}_{1}(R[G])\rightarrow 
\widehat{\Kr}_{1}(R[G])\rightarrow \widehat{\mathrm{
Det}} ( \Kr_{1}(R[G]) ) \rightarrow 1.
\end{equation}
\end{prop}

The most important tool in our study of $\Kr_{1}$ for $p$-adic group rings is
the group logarithm (defined in \S \ref{3.1}). We denote by $G^{\rm ab}$ the abelianization of $G$. 
Also denote by $I(R[G])$ the kernel of the augmentation map 
$R[G]\rightarrow R$ and by $\mathcal{A}(R[G])$  the kernel of the natural 
map from $R[G]$ to $R[G^{\rm ab}] $.
 The following two results are both
proved using this technique:

\begin{thm}\label{thm5}
Let $G$ be a $p$-group. Let $C_{G}$ denote the set of conjugacy classes of $G$ and
let $\phi :R[G] \rightarrow R [ C_{G} ] $ be the $R$-linear map induced by sending each element of $G$ to its conjugacy class. 
The images $\phi(\mathcal{A} (R[G]))$, $\phi(I (R[G]))$ are finite free $R$-modules with
$\phi(\mathcal{A} (R[G]))\subset \phi(I (R[G]))\subset R[C_G]$.

a) There 
is an exact sequence 
\begin{equation*}
1\rightarrow \phi  ( \mathcal{A} (R[G]))
\xrightarrow{\ \mu\ } \mathrm{Det}(R[G]^{\times } )
\rightarrow R [ G^{\rm ab} ]^{\times }\rightarrow 1
\end{equation*}
where the right-hand map is induced by passage to quotient and the left-hand
map $\mu$ is induced by the inverse of the group logarithm map described in
Section \ref{3}.

b) There is an exact sequence
\begin{equation*}
\mathrm{Det} ( R [ G ] ^{\times } ) \xrightarrow{\ {\mu}' \  }
 \phi  ( I(R [ G ] ) ) \xrightarrow { \ } \frac{R}{
 ( 1-F ) R}\otimes _{\mathbf{Z}_{p}}G^{\rm ab}\rightarrow 0.
\end{equation*}
where $\mu'$ is induced by the group logarithm map. 
The composition $\mu'\cdot \mu$ is the inclusion 
$\phi(\mathcal{A} (R[G]))\hookrightarrow \phi(I (R[G]))$.

\end{thm}

For a number of our results which concern an arbitrary finite group $G$ we
shall need to impose the extra condition that $R$ is a Noetherian, normal
ring.

\begin{thm}\label{thm6}
Let $S$ be an integral domain which contains $R$ and which also satisfies
the standing Hypothesis. (But note that we do not need to suppose that ${F_{S}}_{|R}=F_{R}$.) Suppose further that $S$ and $R$ are Noetherian normal rings;
that $S$ is a finitely generated free $R$-module; that the field of
fractions of $S$ is a Galois extension of $N$ with Galois group $\Delta$
which acts on $S$ with the property that $S^{\Delta }=R$; and that the trace
map ${\rm Tr}:S\rightarrow R $ maps $S$ onto $R$. Let $G$ be an arbitrary finite
group and let $\Delta$ act coefficientwise on $\mathrm{Det} ( S [ G]^{\times })$; then 
\begin{equation*}
\mathrm{Det} ( \GL(S[G]))^{\,\Delta }=
\mathrm{Det}( \GL(R[G])) .
\end{equation*}
\end{thm}

This completes our brief summary of our results pertaining to the
determinantal image of $\Kr_{1}(R[G])$. Next we
give a flavour of some forthcoming arithmetic applications of the above
results. Consider a finite group $G$, a normal two dimensional scheme $Y$
 which is projective and flat over $\mathrm{Spec} ( \mathbf{Z} ) $,
and a $G$-cover $\pi :X\rightarrow Y$ with the property that $X$ is a $G$-torsor over $Y$. To such a cover we can associate a projective Euler
characteristic ${\rm R}\Gamma  ( X,\mathcal{O}_{X} ) $ in ${\rm Cl} (\Z[G])$ the projective classgroup of $\Z[G]$. This class is constructed by finding a perfect $\Z[G]$-complex, which is quasi-isomorphic to a Cech complex for $\O_{X}$, which thus calculates the cohomology of $\O_{X}$. One then forms the
alternating sum of the classes of the terms in the perfect complex; see [CE]
for details. We now indicate briefly how the results of this paper may be
used in studying such projective Euler characteristics. For ease of
presentation we cast our discussion in terms of the Zariski topology,
although in practice one will usually work with formal neighbourhoods.

Suppose that $\mathcal{U}=\left\{ U_{i}\right\} $ is an affine cover of $Y$;
we write $U_{ij}=U_{i}\cap U_{j}$, $U_{ijk}=U_{i}\cap U_{j}\cap U_{k}$ and $
U_{i}=\mathrm{Spec} ( R_{i} )$, $U_{ij}=\mathrm{Spec} (
R_{ij})$, $U_{ijk}=\mathrm{Spec} ( R_{ijk} )$. Refining the
cover $\mathcal{U}$ if necessary, we can find local bases $ \{
e_{i} \} $ of $(\pi _{\ast }\mathcal{O}_{X})( U_{i} ) $ over $
R_{i} [ G ]$, so that 
\begin{equation*}
(\pi _{\ast }\mathcal{O}_{X}) ( U_{i} ) =\mathcal{O}_{Y} (
U_{i} )  [ G ] e_{i}=R_{i} [ G ] e_{i};
\end{equation*}%
we can then write $e_{i}=\lambda _{ij}e_{j}$ with $\lambda _{ij}\in R_{ij}
 [ G ] ^{\times }$. The sheaf $\pi _{\ast }\mathcal{O}_{X}$ is said
to have \textit{special} $\mathcal{O}_{Y} [ G ] $-structure if we
can find a cover $\mathcal{U}$ with all the $\mathrm{Det} ( \lambda
_{ij} ) =1$; this then implies that the image of $\lambda _{ij}$ in $ 
\Kr_{1} ( R_{ij} [ G ]  ) $ actually lies in $\SK_{1}(
R_{ij}[ G] )$. We say that  $\pi _{\ast }\mathcal{O}_{X}$
has \textit{elementary} $\mathcal{O}_{Y} [ G ] $-structure if we
can find a cover $\mathcal{U}$ so that each  $\lambda_{ij}$ has trivial
image in $\Kr_{1} ( R_{ij} [ G ]  ) $; this then means that
the $\lambda _{ij}$ (viewed in $\GL( R_{ij}[ G]) $) 
lie in the group of elementary matrices ${\rm E} ( R_{ij} [ G ]
 ) $. The importance of elementary structure is that it allows us to
form a second Chern class as follows. For each triple $i,j,k$ we have the
corresponding Steinberg sequence (see for instance Theorem 5.1 in [M]) 
\begin{equation*}
1\rightarrow \Kr_{2} ( R_{ijk}  G ]  ) \rightarrow {\rm St} (
R_{ijk} [ G ]  ) \rightarrow {\rm E} ( R_{ijk} [ G ]
 ) \rightarrow 1.
\end{equation*}
We choose a section $s:{\rm E}(R_{ijk} [ G ]  ) \rightarrow
{\rm St} ( R_{ijk} [ G ]  ) $; and we then define the 2-cocycle%
\begin{equation*}
z ( i,j,k ) =s ( \lambda _{jk}^{-1} ) s ( \lambda
_{ij}^{-1} ) s ( \lambda _{ik} ) \in \Kr_{2} ( R_{ijk} [
G ]  ) .
\end{equation*}
The $z ( i,j,k ) $ then define a Cech representative for a Zariski
cohomology class $c_2$ in the group ${\rm H}^{2} ( Y,\Kr_{2} ( \O_{Y}[ G]
)) $ which is an equivariant second Chern class for $\pi _{\ast
}\O_{X}$. Because we have put out ourselves in the situation in which the image
of $\lambda _{ij}$ in $\Kr_{1} ( R_{ij} [ G ]  ) $ is
trivial, the natural first Chern class $c_{1}\in {\rm H}^{1} ( Y,\Kr_{1} (
\O_{Y} [ G ]  )  )  $ is trivial. The crucial
importance of this construction is that there is a Riemann-Roch theorem
which enables us to calculate the Euler characteristic ${\rm R}\Gamma ( X,
\mathcal{O}_{X}) $ from the knowledge of the first and second Chern
classes of $\pi _{\ast }\mathcal{O}_{X}$. (See [CPT2].)

We can now explain the role of the results of this paper in this program:
indeed they enable us to show that, under certain circumstances, the $G$-torsor $X$ has special and indeed even elementary structure. For the sake
of simplicity suppose now that the group $G$ is perfect. To each local basis 
$e_{i}$ we associate the resolvend 
\begin{equation*}
r ( e_{i} ) =\sum\limits_{g\in G}\ g\left( e_{i}\right) g^{-1}\in
S_{i}[G]
\end{equation*}
where $S_{i}=(\pi _{\ast }\O_{X})( U_{i}) $ is a finite $R_{i}$-algebra. We note that for $h\in G$ we have 
\begin{equation*}
h( r( e_{i})) =\sum\limits_{g\in G}\ hg (
e_{i} ) g^{-1}=r ( e_{i} ) h.
\end{equation*}
Since $G$ is perfect, we have $\mathrm{Det} ( h ) =1$ and so
\begin{equation*}
h ( \mathrm{Det} ( r ( e_{i})  )  ) =\mathrm{Det}%
( r( e_{i}) h) =\mathrm{Det} ( r (
e_{i} )  ) \cdot \mathrm{Det} ( h ) =\mathrm{Det} (
r ( e_{i} )  ) .
\end{equation*}%
Thus $\mathrm{Det} ( r ( e_{i})) $ is fixed by $G$.
Note that here $G$ acts on $\Kr_{1} ( S_{i} [ G ]  )$, and
hence on $\mathrm{Det} ( S_{i} [ G ] ^{\times } )$, by
its action on the ring $S_i$; that is to say the $G$-action is induced via the
action of $G$ on the coefficients in the group ring $S_{i} [ G ] $.
If we can establish a suitable fixed point theorem, such as in Theorem \ref{thm6},
then we can deduce that $\mathrm{Det} ( r ( e_{i} )  ) =
\mathrm{Det}( \lambda _{i} ) $ for some $\lambda _{i}\in R_{i}
 [ G ] ^{\times }$. If we then replace the bases $ \{
e_{i} \} $ by the bases $ \{ \lambda _{i}^{-1}e_{i} \}$,
then $\mathrm{Det} ( r( e_{i} )  ) =1$. Since it is
easily verified that 
\begin{equation*}
\mathrm{Det} ( r ( e_{i})) =\mathrm{Det} ( \lambda
_{ij} ) \mathrm{Det} ( r ( e_{j} )  )
\end{equation*}
it will follow that $\mathrm{Det} ( \lambda _{ij} ) =1$ and we
shall thereby have exhibited a special structure. At this point we will know
that the image of $\lambda_{ij}$ in $\Kr_{1} ( R_{ij} [ G ]
 ) $ actually lies in $\SK_{1} ( R_{ij} [ G ]  )$.
The approximation theorems on $\SK_{1}$  of group rings, which will be shown
in [CPT1], can then enable us to re-choose our bases to ensure that in fact $\lambda _{ij}$ has trivial image in $\Kr_{1} ( R_{ij} [ G ]
 ) $, and so we finally get our desired elementary structure. 

We conclude this Introduction by briefly describing the structure of the
paper. In Section \ref{2} we recall some standard techniques and results from
K-theory. Section \ref{3} lies at the heart of this paper: here we introduce the
group logarithm, which is a key-tool for our study of $\Kr_{1}$; we then use
this tool to prove Theorems \ref{thm5} and \ref{thm6} for $p$-groups.
The proofs of some of the elementary properties of the group
logarithm in Section \ref{3} are relatively straight forward generalizations of
the corresponding results for the case when $R$ is a non-ramified ring of $p$-adic integers; for such results we either sketch the proof (for the
reader's convenience) or make reference to the existing literature. However, in order to work with the far more general family of coefficient
rings $R$, satisfying the Hypothesis, many aspects of the theory now
become either quite different or at least considerably more subtle. One of
the most striking new features concerns the determination of the image of
the logarithm. To give a flavor of the key-idea here consider the case
when the group $G$ is an abelian $p$-group (see 3.d for full details).
Varying from the above stated convention, we extend the lift of Frobenius $F$
on $ R$ to a lift of Frobenius of the group ring $R [ G ]$, by
stipulating that $F$ is the $p$-th power map on $G$; there is then a natural
extension of $ F$ to the cotangent space (at zero) of $ R [ G ]   $
with the property that $d\cdot F=pF\cdot d.$ In this special situation the group
logarithm is just the composite map $ ( p-F ) \cdot \log $;   thus
the differential of the group logarithm is $d ( p-F ) \cdot \log
=p ( 1-F  ) d\log$ and so the cokernel of the differential of
the group logarithm in $p$ times the cotangent space is given by the group
of coinvariants of  $F$.

The remainder of the paper concerns the use of induction techniques for
reducing the proof of our main results to the corresponding result for $p$-groups. For this we use the method of character action, due in origin to
Lam (see [L]), together with Brauer induction. These methods are introduced
in Section \ref{4}. We then apply these techniques to determinants in Section \ref{5},
where we complete the proofs of Theorems \ref{thm3} and \ref{thm6}.

\bigskip

\section{Generalities}\label{2}
\setcounter{equation}{0}

\subsection{Some $\Kr$-theory}\label{2.1}

In this subsection we gather together a number of standard results that we
shall need in later sections.

Let $S$ denote an arbitrary unitary ring and we write ${\rm M}(S)
=\displaystyle{\lim_{\rightarrow }}({\rm M}_{n}(S))$ where ${\rm M}_n(S)$
stands for $n\times n$-matrices with entries in $S$. Let $I$ denote a two-sided ideal
of $S$ which is contained in the Jacobson radical of $S$; let $\bar{S}
=S/I$. Since all the diagonal entries of a matrix $x\in
1+{\rm M}_{n}(I)$ are units, by left and right multiplying by
elementary matrices we can bring $x$ into diagonal form with unit entries.
Thus $x$ is invertible and so 
\begin{equation*}
1+{\rm M} ( I ) =\GL ( S,I ) \overset{\mathrm{defn}}{=}\ker
 ( \GL ( S ) \rightarrow \GL ( \bar{S} )  ) .
\end{equation*}
Note that $\GL  (S )   $ maps onto $\GL ( \bar{S} ) $,
since given $\bar{x}\in \GL_{n}( \bar{S} ) $ with inverse 
$\bar{y},$ then for lifts $x$, $y\in {\rm M}_{n}(S)$ we have $xy 
\in 1+{\rm M}_{n}(I)$ and so $x$ is invertible.

We have the long exact sequence of K-theory (see for instance page 54 of [M])%
\begin{equation}\label{eq4}
\Kr_{2} ( S ) \xrightarrow{q_{2}} \Kr_{2} (\bar S) \rightarrow \Kr_{1} ( S,I ) \rightarrow \Kr_{1} ( S ) 
\xrightarrow{q_{1}}\Kr_{1}(\bar S)\rightarrow 0
\end{equation}%
with $q_{1}$ surjective since $\GL(S)$ maps onto $\GL(\bar S)$. We let ${\rm E} ( S,I ) $ denote the smallest
normal subgroup of the group of elementary matrices ${\rm E}( S ) $
containing the elementary matrices $e_{ij}( a) $ with $a\in
I$. Then we know (see for instance page 93 in [R])
\begin{equation*}
\Kr_{1} ( S,I ) =\frac{\GL ( S,I ) }{{\rm E} ( S,I ) }.
\end{equation*} 
Under certain circumstances we shall wish to show that $q_{2}:\Kr_{2} (
S ) \rightarrow \Kr_{2} ( \bar S ) $ is surjective (see
for instance 2.3.1). With is in mind we note that by (\ref{eq4}) $q_{2}$ is
surjective if and only if $\Kr_{1} ( S,I ) \rightarrow \Kr_{1} (
S ) $ is injective, and this is the case if and only if
\begin{equation}\label{eq5}
\GL ( S,I ) \cap {\rm E} ( S ) ={\rm E}( S,I ) .
\end{equation}
Recall also from [R] loc. cit. that ${\rm E} ( S,I ) $ is a normal
subgroup of $\GL ( S ) $ and 
\begin{equation}\label{eq6}
\left[ \GL( S) ,{\rm E} ( S,I ) \right] =\left[ {\rm E}(
S ) , {\rm E} ( S,I ) \right] ={\rm E}( S,I ) .
\end{equation}

\begin{lemma}\label{lem7}
Suppose $A=1+\lambda \in 1+{\rm M}_{n} ( I ) =\GL_{n} ( S,I )$. Then 
\begin{equation*}
\left( 
\begin{array}{cc}
A & 0 \\ 
0 & A^{-1}%
\end{array}%
\right) =\left( 
\begin{array}{cc}
1 & \lambda \\ 
0 & 1%
\end{array}%
\right) \left( 
\begin{array}{cc}
1 & 0 \\ 
1 & 1%
\end{array}%
\right) \left( 
\begin{array}{cc}
1 & -A^{-1}\lambda \\ 
0 & 1%
\end{array}%
\right) \left( 
\begin{array}{cc}
1 & 0 \\ 
1 & 1%
\end{array}%
\right) ^{-1}\left( 
\begin{array}{cc}
1 & 0 \\ 
-\lambda & 1%
\end{array}%
\right)
\end{equation*}%
lies in ${\rm E} ( S,I )$.
\end{lemma}

Proof. \ See page 94 in [R].\ \ \ \ \ \ $\square $

\begin{lemma}\label{lem8}
(a) We have the inclusion 
\begin{equation}\label{eq7}
\left[\GL( S) ,1+{\rm M}( I) \right]\subset {\rm E} ( S,I ) ;
\end{equation}

(b) given $g\in \GL ( S,I )$,  there exist $e_{1}$, $ e_{2}\in {\rm E} (
S,I )$, and $x\in 1+I$ such that $g=e_{1}\delta  ( x ) e_{2}$
where $\delta  ( x ) $ is the diagonal matrix with leading term $x$
and with all non-leading diagonal terms equal to $1$; thus in particular 
\begin{equation*}
\GL ( S,I) \subset \left\langle {\rm E}( S,I ) , (
1+I ) \right\rangle
\end{equation*}
where we view $1+I\subset S^{\times }\subset \GL( S,I) $ as
previously.
\end{lemma}

Proof.   We consider terms of the form $x= yry^{-1}r^{-1}$ with 
$r\in \GL_{n} ( S )$ and $y\in 1+{\rm M}_{n} ( I )$; then we
note that 
\begin{equation}\label{eq8}
\left( 
\begin{array}{cc}
yry^{-1}r^{-1 } & 0 \\ 
0 & 1%
\end{array}%
\right) =\left( 
\begin{array}{cc}
y & 0 \\ 
0 & y^{-1}%
\end{array}%
\right) \left( 
\begin{array}{cc}
r & 0 \\ 
0 & 1%
\end{array}%
\right) \left( 
\begin{array}{cc}
y^{-1} & 0 \\ 
0 & y%
\end{array}%
\right) \left( 
\begin{array}{cc}
r^{-1} & 0 \\ 
0 & 1%
\end{array}%
\right) .
\end{equation}%
Now by the previous lemma we know that 
\begin{equation*}
\left( 
\begin{array}{cc}
y & 0 \\ 
0 & y^{-1}%
\end{array}%
\right) \in {\rm E} ( S,I ) .
\end{equation*}%
Thus the first right-hand term in (\ref{eq8}) lies in ${\rm E} ( S,I ) $, and the
product of the last three right-hand terms also lies in ${\rm E} ( S,I ) $, since ${\rm E} ( S,I ) $ is normal in $\GL ( S ) $ hence $x\in {\rm E} ( S,I ) $ as required. 

  To prove (b) consider $g\in \GL ( S,I )$. First we note that,
as ${\rm E} ( S,I )  \subset \GL ( S,I )$ right and left
multiplication of $g$ by elements of ${\rm E} ( S,I )$ yields another
element $\GL ( S,I )$. As explained at the start of this section,
by multiplying $g$ on the right and left by elements of ${\rm E} ( S,I ) $
we may obtain a diagonal matrix all of whose terms are congruent to $1\mod I$; we can multiply this matrix by elements of the form given in
Lemma \ref{lem7} to obtain a diagonal matrix all of whose terms, except the leading
term, are 1.\ \ \ \ \ $\square  $

\subsection{$p$-adic computations}\label{2.2}

In this section we consider various forms of the $p$-adic logarithm. By
definition this logarithm involves certain denominators. In Section 3 we
shall use the lift of Frobenius and the Adams operation $\psi ^{p}$ to
define an integral logarithm, which eliminates these denominators.

Recall that we write $N$ for the field of fractions of $R$ and $C_{G}$
denotes the set of conjugacy classes of the finite group $G$. Let $N
 [ C_{G} ] $ denote the free $N$-vector space on $C_{G}$, and
write $\phi :N [ G ] \rightarrow N [ C_{G} ] $ for the $N$-linear map obtained by mapping each element of $G$ to its conjugacy
class. 

Let $J$ denote a two-sided ideal of $R [ G ] $ with the property
that for some positive integer $m$, $J^{m}\subset pR [ G ]$. Thus
the $p$-adic logarithm yields a map%
\begin{equation*}
\log :1+J\rightarrow N [ G ] ;
\end{equation*}%
moreover, for $n\geq 1$, we have the $p$-adic logarithm map%
\begin{equation*}
\log ( 1+p^{n}R[ G] ) \subset p^{n}R[ G] ;
\end{equation*}
and for $n\geq 2$ we have the (two-sided) inverse map 
\begin{equation*}
\exp :p^{n}R [ G ] \rightarrow 1+p^{n}R [ G ] .
\end{equation*}%
Note that if $ T:R [ G ] \rightarrow \GL_{n} ( R\otimes _{
\mathbf{Z}_{p}}\mathbf{Z}_{p} [ \zeta  ]  ) $ is an integral
representation of $G$ with character $\chi $, then for $x\in 1+J$  we
know that for large enough $N$, $x^{p^{N}}\equiv 1\mod pR[ G
] $ and so $\mathrm{Det}_{\chi }( x) ^{p^{N}}\equiv 1 \mod pR\otimes_{\mathbf{Z}_{p}}\mathbf{Z}_{p}[ \zeta] $;
therefore the usual $p$-adic logarithm $\log ( \mathrm{Det}_{\chi
}( x) )$ is defined in $N\otimes_{\mathbf{Z}_{p}}
\mathbf{Z}_{p}[ \zeta ] $.

\begin{lemma}\label{lem9}
Suppose $x\in 1+J$ and that $\chi $ is any character of $G$ with the
property that $\log \left( \mathrm{Det}_{\chi } ( x ) \right) =0$.
Then $\mathrm{Det}_{\chi }(x)^{p^{M}}=1$ for some positive
integer $M$.
\end{lemma}

Proof. Suppose first that $\chi $ is the character of an integral
representation $T$ as above. Then we know that $\mathrm{Det}_{\chi } (
x ) ^{p^{N+1}}\equiv 1\mod p^{2}R\otimes_{\mathbf{Z}_{p}} \mathbf{Z}_{p}[
\zeta  ] $; however, we have seen that $\log $ is a bijection on $
1+p^{2}R\otimes_{\mathbf{Z}_{p}} \mathbf{Z}_{p}[ \zeta  ] $ and so $\mathrm{Det}%
_{\chi } ( x ) ^{p^{N+1}}=1$. Since monomial representations are
such integral representations, we have now proved the result for the
characters of monomial representations. The result for general characters
then follows by Brauer's induction theorem. $\ \square $

\begin{prop}\label{prop10}
The map $\phi \circ \log $ is a homomorphism and there is a unique
homomorphism $\upsilon $ which makes the following diagram commute
\begin{equation*}
\begin{array}{ccc}
1+J &  \xrightarrow {\mathrm{Det}} & \mathrm{Det} ( 1+J )
\\ 
& \searrow & \downarrow \upsilon \\ 
&  & N[C_{G}]%
\end{array}%
\end{equation*}%
where the diagonal map is $\phi \circ \log$. For $n\geq 2$, $\upsilon $
induces isomorphisms 
\begin{equation*}
\mathrm{Det}( 1+p^{n}R[ G] ) \xrightarrow{\sim} p^{n}R[ C_{G}
],\quad \ \mathrm{Det} ( 1+p^{n}J ) \xrightarrow{\sim} p^{n}\phi  (
J ) .
\end{equation*}
\end{prop}

Proof. For $x\in 1+J$ and for a character $\chi $ of $G$%
\begin{equation*}
\chi \circ \phi \circ \log \left( x\right) =\log \left( \mathrm{Det}_{\chi
}\left( x\right) \right) .
\end{equation*}%
Thus for $y\in 1+J$%
\begin{eqnarray*}
\chi \circ \phi \circ \log \left( xy\right) &=&\log \left( \mathrm{Det}%
_{\chi }\left( xy\right) \right) =\log \left( \mathrm{Det}_{\chi }\left(
x\right) \right) +\log \left( \mathrm{Det}_{\chi }\left( y\right) \right) \\
&=&\chi \circ \phi \circ \log \left( x\right) +\chi \circ \phi \circ \log
\left( y\right)
\end{eqnarray*}%
and, since evaluation on all the irreducible characters of $G$ is injective
on $N[ C_{G}]$, this indeed shows that $\phi \circ \log $ is a
homomorphism.

Suppose now that ${\rm Det}(x)=1$; then $x\in \ker  (
\phi \circ \log  ) $ and this demonstrates the existence of a unique
homomorphism $\upsilon $ which makes the above diagram commute.

To see that the restriction of $\upsilon $ to $\mathrm{Det} ( 1+p^{n}R%
 [ G ]  ) $ is an isomorphism for $n\geq 2$, we show that
the kernels of the restrictions $\phi \circ \log $ and $\mathrm{Det}$
actually coincide. Note that if $x\in \ker  ( \phi \circ \log  ) $,
then for each $\chi $ we know by the above that $\log  ( \mathrm{Det}_{\chi } ( x )  ) =0$, and so by Lemma \ref{lem9} it follows that a $p$-power of 
$\mathrm{Det}_{\chi }(x)$ is equal to one; however,
each 
\begin{equation*}
\mathrm{Det}_{\chi}(x)\equiv 1 \mod p^{n} 
\end{equation*}%
and so, because $n\geq 2$, we can deduce that $\mathrm{Det}_{\chi }(
x) =1$, as required. The argument for $\mathrm{Det}(
1+p^{n}J) $ is entirely similar.\ \ \ \ \ $\square $

\begin{lemma}\label{lem11}
For $n\geq 2$ and for  $x\in ( 1+p^{n}R[ G] ) \cap
\ker ( \mathrm{Det}) $, there exists an element $y$ in $\lbrack R[ G] ^{\times },\ 1+p^{n}R[ G] ]$ so that 
\begin{equation*}
x\equiv y  \mod p^{n+1}R [ G ] .
\end{equation*}
\end{lemma}

Proof. Put $x=1+p^{n}\lambda$. Since $\mathrm{Det} ( x ) =1$ we
know by the above that $\phi \circ \log  ( x ) =0$, and so $\log
( x) =\log( 1+p^{n}\lambda) $ can be written as 
\begin{equation*}
p^{n}\sum\limits_{g,h}\lambda _{g,h}^{\prime }( g-hgh^{-1}) .
\end{equation*}
Since $\log ( x) \equiv p^{n}\lambda \mod p^{n+1}R[ G]$, we have the congruence $\mod p^{n+1}R[ G ] $ 
\begin{equation*}
x=1+p^{n}\lambda \equiv 1+\log  ( x ) \equiv
1+\sum\limits_{g,h}p^{n}\lambda _{g,h}^{\prime } ( g-hgh^{-1} )
\equiv \prod\limits_{g,h} [ h,1-p^{n}\lambda _{g,h}^{\prime }g ]
\end{equation*}%
which gives the required result.\ \ \ \ \ $\square $

\subsection{Some proofs}\label{2.3}

We write $R_{n}$ for the quotient ring $R/p^nR$ and we now apply
the results of \S \ref{2.1} with $S={\rm M} ( R_{n+1} [ G ]  ) $ and $I=p^{n}{\rm M}( R_{n+1}[ G] ) $.

\subsubsection{The surjectivity of $\Kr_{2} ( R_{n+1} [ G ]
 ) \rightarrow \Kr_{2} ( R_{n} [ G ]  ) $.}

The following result does not require $R$ to have a lift of Frobenius.

\begin{prop}\label{prop12}
For $n\geq 2 $ the natural map $\Kr_{2} ( R_{n+1} [ G ]  )
\rightarrow \Kr_{2} ( R_{n} [ G ]  ) $ is surjective.
\end{prop}

Proof. In the following we shall abbreviate $\GL ( R_{n+1} [ G] 
, p^{m}R_{n+1} [ G ]  ) $ to $\GL ( R_{n+1} [ G ]
,p^{m} )$.  By (\ref{eq4}) and (\ref{eq5}) above, it will suffice to show
\begin{equation*}
\GL( R_{n+1}[ G] , p^{n}) \cap [ \GL( R_{n+1}[ G] ), \GL( R_{n+1}[ G] ) ]
\subset {\rm E} ( R_{n+1} [ G ] ,p^{n} ) ,
\end{equation*}
since by the Whitehead Lemma $\left[ \GL( R_{n+1}[ G] )
, \GL ( R_{n+1}[ G] ) \right] ={\rm E}( R_{n+1} [ G
])$. Consider an element $\bar{x}$ of the left-hand
side. We may then choose 
\begin{equation*}
x\in \GL( R[ G] ,p^{n}) \cap [ \GL( R [ G%
 ]  ) , \GL ( R [ G ]  )  ]
\end{equation*}%
with the property that $x\mod p^{n+1}=\bar x$. By Lemma \ref{lem8} (b) we
can find $ e_{1}$, $e_{2}\in {\rm E} ( R [ G ] ,p^{n} ) $ such
that $x=e_{1}de_{2}$ with $d\in  ( 1+p^{n}R [ G ]  ) \cap
\ker  ( \mathrm{Det} )$. So by Lemma \ref{lem11}, $d$ is congruent $\mod p^{n+1}$ to an element of $ [ R [ G ] ^{\times },1+p^{n}R
 [ G ]  ] $ and Lemma \ref{lem8}(a) shows that this element lies in ${\rm E} ( R ( G ) ,p^{n} ) $.\ \ \ \ $\square $

\subsubsection{Proof of Proposition \ref{prop4}}

 We must show that the natural map 
\begin{equation*}
\theta :\mathrm{Det} ( \GL ( R [ G ]  ) )
\rightarrow \widehat{\mathrm{Det}} ( \GL ( R [ G ]  )
 ) =\lim_{\leftarrow }\frac{\mathrm{Det} ( \GL ( R [ G ]
 )  ) }{\mathrm{Det} ( 1+p^{n}{\rm M} ( R [ G ] )
 ) }
\end{equation*}
is an isomorphism.

Let $N^{\prime }$ denote a sufficiently large finite extension of$\ N\ $with
the property that 
\begin{equation*}
\Kr_{1} ( N^{\prime } [ G ] ) =\mathrm{Det} ( \GL (
N^{\prime } [ G])) =\prod N^{\prime \times }.
\end{equation*}

Since $\prod N^{\prime \times }$ is $p$-adically complete, we know that 
\begin{equation*}
\cap_{n}\mathrm{Det}( 1+p^{n}{\rm M}( R [ G ] ))
=\{1\}
\end{equation*}
and so $\theta $ is certainly injective. To show that the image of $\theta $
is surjective, consider a coherent sequence 
\begin{equation*}
\{d_{n}\mathrm{Det} ( 1+p^{n}{\rm M} ( R [ G ] ) )
\}\in \prod\limits_{n}\frac{\mathrm{Det} ( \GL ( R [ G ]
 )  ) }{\mathrm{Det} ( 1+p^{n}{\rm M} ( R [ G ] )
 ) }.
\end{equation*}
It will suffice to show that $\{d_{n}\}$ converges to an element $d\in 
\mathrm{Det} ( \GL ( R [ G]))$. To this end
we consider the sequence $ \{ d_{n}^{\prime } \} = \{
d_{n}d_{2}^{-1} \} $ for $n\geq 2$. We then use Lemma \ref{lem8} (b) to deduce
that for $n\geq 2$ 
\begin{equation*}
\mathrm{Det} ( 1+p^{n}{\rm M} ( R [ G]))=\mathrm{
Det} ( 1+p^{n}R [ G ])
\end{equation*}%
and by Proposition \ref{prop10} the map $\upsilon $ of gives an isomorphism $\mathrm{%
Det} ( 1+p^{n}R [ G ] ) \cong p^{n}R [ C_{G} ]$;
so that we now have isomorphisms 
\begin{equation*}
\delta _{n}:\mathrm{Det} ( 1+p^{n}{\rm M} ( R [ G])
 ) \cong p^{n}R [ C_{G} ] .
\end{equation*}%
It therefore follows that the sequence $ \{ \delta _{2} (
d_{n}^{\prime } )  \} $ is a Cauchy sequence in $p^{2}R [ C_{G}
 ]$, which converges to an element $x\in  p^{2}R [ C_{G} ] $. The desired element $d$ is then given by $d_{2}\delta _{2}^{-1} (
x )$.

It now remains to construct the exact sequence (\ref{eq3}). First we note that by
definition $\widehat{\SK}_{1}( R[ G] ) $ is the inverse
limit of the images of $\SK_{1}( R[ G] ) $ in the
groups $\widehat{\Kr}_{1} ( R_{n} [ G ] ) ,$ which we
denote by $\SK_{1}(R[G])_{n}$. By the localisation
sequence (\ref{eq4}) we know that 
\begin{equation*}
K_{1} ( R_{n} [ G ]  ) =\frac{\Kr_{1} ( R[G])}{\Kr_{1} ( R[G] , p^{n} ) }
\end{equation*}
and by Lemma \ref{lem8} (b) this is given by the image of $\GL ( R [ G ]
 ) $ modulo the image of $1+p^{n}{\rm M} ( R [ G ]  ) $.
Thus if we write $\SK_{1}( R[G] ,p^{n} ) $ for the
intersection of $\Kr_{1} ( R[ G ] ,p^{n} ) $ and $\SK_{1}(R[G]) $ then we have a commutative diagram
with exact rows:
\begin{equation*}
\begin{array}{ccccccccc}
1 & \rightarrow & \SK_{1}( R[ G] ,p^{n}) & \rightarrow
& \Kr_{1} ( R[G] ,p^{n} ) & \rightarrow & \mathrm{Det}%
 ( 1+p^{n}M ( R [ G ]  )  ) & \rightarrow & 1 \\ 
&  & \downarrow &  & \downarrow &  & \downarrow &  &  \\ 
1 & \rightarrow & \SK_{1} ( R [ G]) & \rightarrow & 
\Kr_{1} ( R [ G])) & \rightarrow & \mathrm{Det} (
\GL ( R [ G ]  )  ) & \rightarrow & 1 \\ 
&  & \downarrow &  & \downarrow &  & \downarrow &  &  \\ 
1 & \rightarrow & \SK_{1} ( R [ G ]  ) _{n} & \rightarrow & 
\Kr_{1} ( R_{n} [ G ] ) & \rightarrow & \mathrm{Det} (
\GL ( R [ G ]  )  ) _{n} & \rightarrow & 1%
\end{array}%
\end{equation*}%
where we write $\mathrm{Det} ( \GL ( R [ G ]  )  )
_{n}$ for quotient of $\mathrm{Det} ( \GL(R[G]))$ by $\mathrm{Det} ( 1+p^{n}{\rm M} ( R [ G ]  )
 )$. We then take the inverse limit of the bottom exact row; noting
that the the inverse limit of the $\SK_{1}( R[G])_{n} $ affords a surjective system, we obtain the desired exact sequence
(\ref{eq3}). \ $\square $
\bigskip

\section{The group logarithm and determinants for $p$-groups}\label{3}
\setcounter{equation}{0}

\subsection{The group logarithm.}\label{3.1}

In this subsection we introduce and develop the group logarithm which is a
fundamental tool for studying the determinants of units of $p$-adic group
rings. Many of the proofs of the results we give in this subsection are
comparatively routine extensions of the corresponding results for group
rings whose coefficients are rings of integers of $p$-adic fields (see for
instance [T], [F], [CR2] and [O2]); when this is the case we shall confine
ourselves to giving a brief idea of the proof (for the reader's convenience)
and provide a  precise reference for the result over a $p$-adic ring of
integers. \medskip

Recall that $F$ denotes a lift of Frobenius on $R$. We let $G$ denote a
finite $p$-group, and $I(R[G])$ denotes the
augmentation ideal of the group ring $R[G]$. Recall also that
for a left Artinian ring, the Jacobson radical is a nilpotent ideal (see for
instance 5.15 in [CR1]); so, because $I(\mathbf{F}_{p}[G])$ is
the Jacobson radical of the Artinian ring $\mathbf{F}_{p}[G]$,
it follows that we can find a positive integer $N$ such that $I(\mathbf{F}_{p}[G])^{N}=0$. Recall $R$ is torsion free; so $I(R [
G ] )=R\cdot I(\mathbf{Z}_{p} [ G ] )$, and it therefore follows
that 
\begin{equation}\label{eq9}
I(R [ G ] )^{N}\subset pR [ G ] .
\end{equation}

\begin{prop}\label{prop13}
(a) Suppose $\mathrm{Det}(x) \in \mathrm{Det} ( \GL_{n} (
R[G])) $ has the property that $$\mathrm{Det} (
x )  ( {1}_{G} ) =1, $$ where ${1}_{G}$ denotes
the trivial character of $G$. Then 
$$
\mathrm{Det}(x) \in \mathrm{%
Det} ( \GL_{n}(R[G]) ,I(R[G])))
 =\mathrm{Det} ( 1+I ( R [ G]  )
 ).$$

(b) We have $\mathrm{Det} ( \GL_{n} ( R [ G])) =%
\mathrm{Det} ( R [ G ]^{\times } )$ and this proves Theorem 
\ref{thm3} (a).
\end{prop}

Proof. Let $\varepsilon :R[G] \rightarrow R$ denote the
augmentation map, and for each $n\geq 1$, let $T_{n}:\GL_{n}( R[ G%
] ) \rightarrow \GL_{n}( R) $ be the representation
induced by $\varepsilon$. Then we have the split exact sequence
\begin{equation*}
1\rightarrow \GL_{n}( R[G] , I ( R[G])) \subset \GL_{n} ( R [ G ]  ) \leftrightarrows \GL_{n} ( R ) \rightarrow 1.
\end{equation*}%
By hypothesis $\mathrm{Det}( T_{n}( x)) =1$ and so,
if we put $y=x\cdot T_{n} ( x^{-1} )$, then $T_{n}(y)
=1_{n}$, and therefore $y\in \GL_{n} ( R [ G ] ,I ( R [
G ]  ) ) $ and (a) follows since $\mathrm{Det} (
x ) =\mathrm{Det} ( y ) $. To conclude the proof of (a), we
use the fact that, since $I ( R [ G]) $ is contained in
the Jacobson radical of $R[G]$, by Lemma \ref{lem8} (b) we know that 
\begin{equation*}
\mathrm{Det} ( \GL_{n} ( R [ G ] , I ( R [ G ]
 )  )  ) =\mathrm{Det} ( 1+I ( R [ G ]
 )  ) .
\end{equation*}
Part (b) follows from the exact sequences:
\begin{equation*}
\begin{array}{ccccccccc}
1 & \rightarrow & \mathrm{Det} ( \GL_{n} ( R [ G ] ,I ( R[ G]))) & \rightarrow & \mathrm{Det} (
\GL_{n} ( R [ G])) & \rightarrow & \mathrm{Det}%
 ( \GL_{n}(R)) & \rightarrow & 1 \\ 
&  & \uparrow &  & \uparrow &  & \ \ \uparrow \det^{-1} &  &  \\ 
1 & \rightarrow & \mathrm{Det} ( 1+I ( R [ G])
 ) & \rightarrow & \mathrm{Det} ( R [ G ] ^{\times } )
& \rightarrow &R^{\times }\ \ \  & \rightarrow & 1%
\end{array}%
\end{equation*}%
The right-hand vertical map is an isomorphism, its inverse being given by
the determinant for matrices over the commutative ring $R;$ we have just
seen that the left-hand vertical map is an isomorphism; so finally we deduce
that the middle vertical map is also an isomorphism.\ \ \ \ \ $\square
\medskip $

Define the $F$-semi-linear map $\Psi :R[G] \rightarrow R[G]$ by the rule that 
$\Psi  ( rg ) =F ( r ) g^{p}$; as previously $\phi :R [ G ] \rightarrow R [ C_{G} ] $
denotes the $R$-linear map given by sending each group element to its
conjugacy class. Write $\bar{\Psi }:R [ C_{G} ] \rightarrow R [ C_{G} ] $ for the $F$-semi-linear map induced by $\Psi $.

We now define the group logarithm 
$$
\mathcal{L}:1+I(R [ G ] )\rightarrow
N [ C_{G} ]
$$
by the rule that for $x\in I(R[G])$
\begin{eqnarray}\label{eq10}
\mathcal{L} ( 1-x ) &=& ( p-\bar{\Psi } ) \phi  (
 ( \log  ( 1-x )  )  ) =\phi (  ( p-\Psi
 ) ( \log  ( 1-x )  )  ) \\
&=&-\phi  \left( \sum_{n\geq 1}\frac{px^{n}}{n}-\sum_{n\geq 1}\frac{\Psi
 ( x^{n} ) }{n} \right)  \notag
\end{eqnarray}
which is seen to converge to an element of $N [ C_{G} ] $ by (\ref{eq9}).
Note that the map $\mathcal{L}=\mathcal{L}_{F}$ depends on the chosen lift
of Frobenius $F$; we shall therefore have to be particularly careful when
using such group logarithms for different coefficient rings which may well
have different lifts of Frobenius.

\begin{lemma}\label{lem14}
For a character $\chi $ of $G$ and for $x\in I( R[ G]) $
\begin{equation*}
\chi ( \phi ( ( \log ( 1+x) ) )
) =\log ( \mathrm{Det}( 1+x)( \chi))
\end{equation*}%
and 
\begin{equation*}
\chi ( \mathcal{L}( 1+x) ) =\log \, [ \mathrm{Det}%
( 1+x) ( p\chi ) \mathrm{Det}( 1+F( x)
) ( -\psi ^{p}\chi ) ]
\end{equation*}%
where $\psi^{p}$ denotes the $p$th Adams operation on virtual characters of 
$G$, which is defined by the rule $(\psi ^{p}\chi) ( g) =\chi
( g^{p})$.
\end{lemma}

Proof. The key point here is that the logarithm of the determinant is the
trace of the logarithm, as seen in the proof of Proposition \ref{prop10}; for full
details see Proposition 1.3 on page 53 of [T].\ \ \ $\square $

\begin{cor}\label{cor15}
If $\mathrm{Det} ( 1+x ) =1$, then $\phi ( \log  (
1+x)) =0,$ and so $\mathcal{L}( 1+x ) =1$. Thus there
is a unique map 
$$
\nu :\mathrm{Det}\left( 1+I(R [ G\right]))
\rightarrow N[ C_{G}]
$$
such that $\mathcal{L}=\nu \circ \mathrm{Det}$.
\end{cor}

By a slight abuse of language, we will sometimes also refer to $\nu$ as
the group logarithm.
\smallskip

Proof. For such $x$ we know from the proof of Proposition \ref{prop10} that 
we have $$\phi
\left( \left( \log \left( 1+x\right) \right) \right) =0;
$$
 hence
\begin{equation*}
\mathcal{L}\left( 1+x\right) =\left( p-\overline{\Psi }\right) \phi \circ
\log \left( 1+x\right) =0.\ \ \ \square 
\end{equation*}
\medskip

We now show that the group logarithm is an integral logarithm.

\begin{thm}\label{thm16}
(a) We have the inclusion%
\begin{equation*}
\mathcal{L} \, ( 1+I(R [ G ] ) ) \subset p\, \phi  ( I (
R [ G ]  )  ).
\end{equation*}

(b) If $G$ is abelian, then $ \mathcal{L} ( 1+I^{2}(R [ G ]
) ) \subset p\, I^{2}(R [ G ] )$.
\end{thm}

We only provide a sketch of the proof. (See Theorem 54.5 in [CR2] and Theorem 16 in [F] for details.)
We rewrite the right-hand side of (\ref{eq10}) as 
\begin{equation*}
-\phi \left( \left(\sum_{n\geq 1,p\nmid n}\frac{px^{n}}{n}+\sum_{n\geq 1}\frac{
px^{np}}{np}\right)-\sum_{n\geq 1}\frac{\Psi ( x^{n}) }{n}\right)
\end{equation*}%
and so the theorem will follow from the following: for any $x\in R[G]$ and $m\geq 0$
\begin{equation*}
\phi  \left( x^{p^{m+1}}-\Psi  ( x^{p^{m}} ) \right) \in p^{m+1}R[ C_{G}]
\end{equation*}
and this follows from the \textquotedblleft non-commutative binomial
theorem\textquotedblright\ as given in Lemma 5.2 on page 80 of [F]. \ This
shows that $\mathcal{L}\left( 1+I(R [ G ] )\right) \subset p\, R [
C_{G} ]$. From Lemma \ref{lem14} we know that $\mathcal{L} ( 1+I(R [ G%
 ] ) ) $ vanishes under evaluation by the trivial character ${1}_{G}$ and so $\mathcal{L} ( 1+I(R [ G ] ) )
\subset p\, \phi  ( I ( R [ G ] ))$.

The proof of (b) is similar but easier. Suppose now that $G$ is abelian;
then $\Psi $ is a ring homomorphism and $\Psi ( I ( R [ G ])) \subset I ( R [ G ])$. Consider $
\left\{ a_{h},b_{h}\right\} \in I( R [ G ]  ) $ for $
h=1,\cdots,  n$. Then by the binomial theorem for each $h$ 
\begin{equation*}
a_{h}^{p}\equiv \Psi  ( a_{h} ) \mod p\, I ( R [ G ]
 ) ,\ \ b_{h}^{p}\equiv \Psi \left( b_{h}\right)  \mod p\, I ( R%
 [ G ]  )
\end{equation*}%
and so 
\begin{equation*}
\left( \sum\limits_{h=1}^{n}a_{h}b_{h}\right)^{p}\equiv \Psi \left(
\sum\limits_{h=1}^{n}a_{h}b_{h}\right) \mod p\, I^{2} ( R [ G]).
\end{equation*}%
Hence for $m>0$%
\begin{equation*}
\left( \sum\limits_{h=1}^{n}a_{h}b_{h}\right)^{p^{m+1}}\equiv \Psi \left(
\sum\limits_{h=1}^{n}a_{h}b_{h}\right) ^{p^{m}}\mod p^{m+1}I^{2} (R[G])
\end{equation*}
and so we have shown that for any $x\in I^{2} ( R[G])$ 
\begin{equation*}
x^{p^{m+1}}\equiv \Psi \left( x^{p^{m}}\right) \mod p^{m+1}I^{2}(R[G]).
\end{equation*}%
Then, arguing as in part (a), we see that $\mathcal{L} ( 1+I^{2} ( R[G])) \subset I^{2} ( R[G])$.\ \ \   $\square $

\bigskip

To gain a more precise idea of the image of the group logarithm, we set 
$$
\mathcal{A} ( R[G]) ={\rm Ker} ( R[G])
\rightarrow R [ G^{\rm ab}])
$$
 and we note that obviously $
\mathcal{A} ( R[G])) \subset I ( R [ G ])$. Our next aim is to show that

\begin{thm}\label{thm17}
We have the equality
\begin{equation*}
\mathcal{L} ( 1+\mathcal{A}(R[G])) =p\, \phi  ( 
\mathcal{A}(R[G])) .
\end{equation*}
\end{thm}

Let $c$ denote a central element of $G$ of order $p$ and set $\overline G=G/\langle c\rangle$.

\begin{lemma}\label{lem18}
(a) We have $$
\left( 1-c\right) ^{p}\equiv -p\left( 1-c\right) \mod p\left(
1-c\right) ^{2}.
$$

(b) For $\xi \in R[G]$, we have 
\begin{eqnarray*}
\mathcal{L} ( 1+ ( 1-c ) \xi  ) &\equiv &p\, \phi  (
 ( \xi -\xi ^{p} )  ( 1-c )  ) \mod p\, \phi  (
 ( 1-c ) ^{2}R [ G ]  ) \\
&\equiv &p\, ( \phi  ( \xi  ) -\overline{\Psi }\circ \phi  (
\xi  )  ) ( 1-c ) \mod p\phi  (  (1-c ) ^{2}R [ G ]  ).
\end{eqnarray*}
\end{lemma}

Proof. The proof when $p=2$ is trivial, so now we assume that $p$ is odd.\
By the Binomial theorem 
\begin{equation*}
1=c^{p}=\left( \left( c-1\right) +1\right) ^{p}\equiv \left( c-1\right)
^{p}+p\left( c-1\right) +1\mod p\, \left( 1-c\right) ^{2}.
\end{equation*}
To prove the second part, we note that using the first part of the lemma we
have the following congruences $\mod p\, ( 1-c ) ^{2}$ 
\begin{eqnarray*}
\mathcal{L} ( 1+ ( 1-c ) \xi  ) &\equiv &p\,\phi  (
 ( 1-c ) \xi  ) + ( 1-c ) ^{p}\phi  ( \xi
^{p} ) \\
&\equiv &p\,( 1-c ) \phi  ( \xi  ) -p ( 1-c ) \phi
 ( \xi ^{p} ) ;
\end{eqnarray*}%
the result then follows since $\phi  ( \xi ^{p} ) \equiv \bar{
\Psi }\circ \phi  ( \xi  ) \mod p$, as we have also seen above.\ \ \ $\square $ 
\medskip

\begin{lemma} \label{lem19}
We have

(a) 
 \ \ $\Psi (  ( 1-c ) R [ G ]  ) =0$.
 
(b)
 \ \ $
\mathcal{L} ( 1+ ( 1-c ) R [ G ]  ) \subset p\phi
 (  ( 1-c ) R [ G ]  ) .
$
\end{lemma}

Proof. As $\Psi  ( 1-c ) =1-c^{p}=0$ and as $c$ is central, it
follows that $\Psi  (  ( 1-c ) R[ G ]  ) =0.$
So, using the above fact that $ ( 1-c ) ^{p}\in p ( 1-c ) R%
 [ G ] $,\ we see that for $x\in R [ G ] $%
\begin{eqnarray}\label{eq11}
\mathcal{L} ( 1- ( 1-c ) x ) &=&-p\, \phi \left( \sum_{n\geq
1}\frac{ ( 1-c ) ^{n}x^{n}}{n} \right) \\
&=&-\phi  ( p ( 1-c ) x+y )  \notag
\end{eqnarray}%
with $y\in p\, ( 1-c ) x^{2}R [ G ]$.\ \ $ \square $

\begin{lemma} \label{lem20}

(a) We have \ $
\mathcal{L}( 1+( 1-c) I( R[G])) =p\phi (( 1-c) I ( R[G])) .
$

(b) For $n\geq 2$, we have
$
\mathcal{L}( 1+( 1-c) ^{n}R[ G]) =p\, \phi
( ( 1-c) ^{n}R[G]).
$
\end{lemma}

Proof. Both results follow by successive approximation, on using the last
of the expressions in (\ref{eq11}) and taking $x$ to lie in $I^{m}(R[G]) $ for (a) 
(resp. $ ( 1-c) ^{n-1}I^{m}( R[G]) $ \ for (b)) with increasing $ m$. \ \ \ \ \ $ \square   $

Note that the central cyclic group $\langle c\rangle $ acts on
the set of conjugacy classes $C_{G}$ by multiplication.

\begin{lemma}\label{lem21}
\begin{eqnarray}\label{eq12-13}
\left( 1-c\right) R[ C_{G}] &=&\ker \left( R[ C_{G}]
\rightarrow R[ C_{\overline{G}}] \right) \\
&=&\ker \left( \phi ( I( R[ G]))
\rightarrow \bar{\phi }( I( R[ \overline{G}]))\right)
\end{eqnarray}
\begin{equation}\label{eq14}
\ p^{n}R[ C_{G}] \cap  ( 1-c ) N [ C_{G} ]
=p^{n}( 1-c) R[ C_{G} ] =\phi ( p^{n}(
1-c) R[ G] ) .
\end{equation}
Furthermore, if $\ c\in \left[ G,G\right] $, then 
\begin{equation}\label{eq15}
\left( 1-c\right) R\left[ C_{G}\right] =\ker \left( \phi ( I^{2}( R[ G])) \rightarrow \overline{\phi }(
I^{2}( R[ \overline{G}] )) \right)
\end{equation}
\end{lemma}

Proof. To prove the first equality we note that, since $ \langle
c \rangle $ acts on $C_{G},$ $R[ C_{G}] $ is a permutation $R\langle c\rangle $-module, and so $R[ C_{G}]  $ is
isomorphic to: a sum of copies of $R$ (coming from the conjugacy classes of
elements $g\in G$ with the property that the elements $\left\{
c^{i}g\right\} _{i=0}^{p-1}$ are all conjugate); and a sum of copies of $%
R\langle c\rangle \ $(coming from the conjugacy classes of
elements $g\in G$ with the property that the conjugacy classes of the
elements $\left\{ c^{i}g\right\} _{i=0}^{p-1}$ are all distinct). The
equality then follows by considering these two cases. The third equality
then also follows since we have now shown that $ ( 1-c) R[C_{G}] $ is an $R$-direct summand of $R[ C_{G}]$.
The second equality follows from the fact that $\left( 1-c\right) R[
C_{G}] \subset \phi ( I( R[ G]))$.
Similarly to show the fourth equality, it will suffice to show that if $c\in %
\left[ G,G\right] $ then 
\begin{equation*}
\left( 1-c\right) R[ G] \subset \ker \left( I^{2}( R[ G] ) \rightarrow I^{2}( R[ \overline{G}] )
\right) .
\end{equation*}
This follows at once from the congruence $\left( 1-g\right) +\left(
1-h\right) \equiv \left( 1-gh\right) \mod I^{2}( R[ G]) $ for $g$, $h\in G$.\ \ \ \ $\square $

Note that, if we wish, we may impose the further condition that $c$ be a
commutator. Indeed, we know that the last term in the lower central series
of $G$ is generated by central commutators; we choose a non-trivial such
element $d=\left[ g,h\right] $ and we suppose that $d$ has order $p^{m}$;
then $dh=ghg^{-1}$, and so$\ d^{p^{m-1}}h^{p^{m-1}}=gh^{p^{m-1}}g^{-1}$,
and hence $c=\left[ g,h^{p^{m}}\right] $ is a central commutator of order $
p $. As $c$ has trivial image in $G^{\rm ab},\ $note that: 
\begin{equation*}
\left( 1-c\right) R[ G] \subset \mathcal{A}(R[ G] ).
\end{equation*}

\begin{lemma}\label{lem22}
If $c$ is a central commutator of order $p,$ then%
\begin{equation*}
\mathcal{L}\left( 1+\left( 1-c\right) R[ G] \right) =p\,\phi \left(
\left( 1-c\right) R[ G] \right) =p\,\phi \left( \left( 1-c\right)
I\left( R[ G] \right) \right) .
\end{equation*}
\end{lemma}

Proof. Put $c=\left[ \gamma ,\delta \right]$. As $R[G]
=R\delta +I(R[G])$, 
\begin{eqnarray*}
\left( 1-c\right) R[ G] &=&R\left( 1-c\right) \delta +\left(
1-c\right) I( R[ G] ) \\
&=&R\left( 1-\gamma \delta \gamma ^{-1}\delta ^{-1}\right) \delta +\left(
1-c\right) I( R[ G] ) \\
&=&R\left( \delta -\gamma \delta \gamma ^{-1}\right) +\left( 1-c\right)
I( R[ G])
\end{eqnarray*}%
and the second equality follows because $\phi ( \delta -\gamma \delta
\gamma ^{-1}) =0$. The first equality now follows from the above
together with Lemmas \ref{lem19} and \ref{lem20} (a). \ \ \ $\square$
\medskip

\textit{Sketch of the proof of Theorem \ref{thm17}.}\ (For additional details see the proof of
Proposition 1.11 in Chapter 6 of [T].)  We shall prove the theorem by
induction on the order of $G.$ Note that the theorem is trivial if $G$ is
abelian; this then starts the induction. Put $\overline{G}=G/\langle
c\rangle $ with $c$ a central commutator of order $p$ as above. Then
by the induction hypothesis we know that the theorem is true for $\overline{G}$.

The proof of the theorem proceeds in two steps.
\smallskip

{\sl Step 1.} We show $\mathcal{L}( 1+\mathcal{A}(R[ G] ))
\supseteq p\,\phi ( \mathcal{A}(R[ G] ))$.

Given $a\in \mathcal{A} ( R [ G ]  )$, since $\mathcal{A}
 ( R [ G ]  ) $ maps onto $\mathcal{A} ( R [ 
\overline{G} ]  ) $, using the theorem for $\overline{G}$, we can
find $x\in \mathcal{A} ( R [ G ]  ) $ such that $\mathcal{L%
} ( 1+\overline{x} ) =p\,\phi  ( \overline{a} )$. Let $
y=p\phi  ( a ) -\mathcal{L} ( 1+x ) $; by Theorem \ref{thm16} we
know that $y\in p\,R [ C_{G} ]$, and so by (\ref{eq12-13}) 
\begin{equation*}
y=p\, \phi (a) -\mathcal{L} ( 1+x ) \in p\,\ker ( R%
 [ C_{G} ] \rightarrow R [ C_{\overline{G}} ]  )
=p\, ( 1-c ) R [ C_{G} ] .
\end{equation*}%
Now by Lemma \ref{lem22}, we can find $b\in  ( 1-c ) R [ G ]
\subset \mathcal{A}(R [ G ] )$ such that $\mathcal{L} (
1+b ) =y;$ therefore 
\begin{equation*}
\mathcal{L} ( ( 1+x )  ( 1+b )  ) =p\,\phi  (
a ) .
\end{equation*}
\smallskip

{\sl Step 2.} We show $\mathcal{L}\left( 1+\mathcal{A}(R\left[ G\right] )\right)
\subset p\phi \left( \mathcal{A}(R\left[ G\right] )\right) .$

Given $a\in \mathcal{A}(R[G])$, using the theorem
for $\overline{G}$, we can find $b\in \mathcal{A}(R[G]) $ such that $\mathcal{L} ( 1+\overline{a} ) =p\, \phi  ( 
\overline{b} ) .$ By Theorem \ref{thm16} we know that $\mathcal{L} (
1+a ) \in p\,R [ C_{G} ] $ and so using (\ref{eq12-13}) we get 
\begin{equation*}
\mathcal{L} ( 1+a ) -p\phi  ( b ) \in p\ker  ( R [
C_{G} ] \rightarrow R [ C_{\overline{G}} ] ) =p (
1-c ) R [ C_{G} ] \subset p\, \phi  ( \mathcal{A}(R [ G%
 ] ) ) . \ \ \square 
\end{equation*}

\begin{prop}\label{prop23}
We have the equality%
\begin{equation*}
\mathcal{L}\left( 1+I^{2} ( R [ G ]  ) \right) =p\phi\left( I^{2} ( R [ G ]  ) \right) .
\end{equation*}
\end{prop}

Proof. As a first step we show that $\mathcal{L}\left( 1+I^{2} ( R [
G ]  ) \right) \subset p\phi \left( I^{2} ( R [ G ]
 ) \right)$. For this we first note that $\mathcal{A} ( R [ G]  ) =\ker \left( R [ G ] \rightarrow R [ G^{\rm ab}%
 ] \right) $ is generated by terms of the form $g ( \left[ h,k%
\right] -1 ) $ for $g$, $h$, $k\in G$ and so as in the proof of Lemma \ref{lem21}
above we see that $\mathcal{A}(R[G]) \subset
I^{2} ( R [ G])$. Now consider $x\in I^{2}(R[G])$. By Theorem \ref{thm16} we can find $y\in I^{2}(R[G])$ such that 
\begin{equation*}
\mathcal{L} ( 1+x ) -p\phi  ( y ) \in \ker \left( pR [
C_{G} ] \rightarrow pR [ G^{\rm ab} ]  \right) =p\phi  \left( 
\mathcal{A}(R[G])\right) \subset p\phi \left(
I^{2}(R[G])\right)
\end{equation*}
and so $\mathcal{L} ( 1+x ) \in p\phi  \left( I^{2}(R[G]\right)$.
\smallskip

We are now in a position to prove the theorem by induction on the group
order. Suppose first that $G$ is abelian and let $c\in G$ have order $p$.
Then we have a commutative diagram with exact rows, where the vertical
arrows are induced by the group logarithm:
\begin{equation*}
\begin{array}{ccccccccc}
1 & \rightarrow & 1+\left( 1-c\right) I(R[G])  & 
\rightarrow & 1+I^{2}(R[G])  & \rightarrow & 
1+I^{2}(R[\overline {G}])  & \rightarrow & 1 \\ 
&  & \downarrow &  & \downarrow &  & \downarrow &  &  \\ 
0 & \rightarrow & p\left( 1-c\right) I(R[G]) & 
\rightarrow & p\, I^{2}(R[G]) & \rightarrow & 
p\, I^{2}(R[\overline {G}]) & \rightarrow & 0.%
\end{array}%
\end{equation*}%
The right-hand map is surjective by the induction hypothesis; the left-hand
vertical map is surjective by Lemma \ref{lem20} (a); and so the central vertical map
is surjective.

Suppose now that $G$ is non-abelian. Then, as above, we let $c$ denote a
central element of $G$ of order $p$, which is a commutator, and put $
\overline{G}=G/ \langle c \rangle$. Then we have a commutative
diagram where the vertical arrows are induced by the group logarithm:
\begin{equation*}
\begin{array}{ccccccccc}
1 & \rightarrow & 1+\left( 1-c\right) R[G]& \rightarrow & 
1+I^{2} ( R[G]) & \rightarrow & 1+I^{2} (R[\overline{G}]) & \rightarrow & 1 \\ 
&  & \downarrow &  & \downarrow &  & \downarrow &  &  \\ 
0 & \rightarrow & p\phi \left( (1-c) R[G]\right) & 
\rightarrow & p\phi \left( I^{2}(R[G]) \right) & 
\rightarrow & p\overline{\phi }\left( I^{2}(R[\overline {G}]) \right) & \rightarrow & 0.%
\end{array}%
\end{equation*}%
The bottom row is exact by Lemma \ref{lem21}, and indeed the proof of this result
shows that the top row is exact. We then conclude as previously, but now
appealing to Lemma \ref{lem22}, to show that the left-hand vertical map is
surjective.\ \ \ \ $\square  $
\medskip

 {\sl When $c$ is not a commutator.} In the above we needed to
choose $c$ to be a commutator. We can always do this, and indeed, as we have
seen, this case provides the crucial inductive step for building up the
group logarithm. However, when we come to our treatment of $\SK_{1}$ in
[CPT1], we shall also need results for the case where $c$ is central of
order $p$ but is not a commutator; we now consider this case.

\begin{lemma}\label{lem24}
If $c$ cannot be written as a commutator in $G$, then multiplication by $c$
permutes the elements of $C_{G}$ without fixed points, and the kernel of
multiplication by $1-c$ on $ R [ C_{G} ] $ is generated by elements
of the form $\sum_{i=0}^{p-1}c^{i}\phi ( g) $ for $g\in G$. 
\end{lemma}

Proof. Suppose to the contrary that $c\phi(h) =\phi(h)$ for $h\in G$. Then $ch=ghg^{-1}$ for some element $g\in G$ and so 
$c=\left[ g,h\right]$, which is the desired contradiction. As previously,
we then observe that a permutation $R \langle c \rangle $-module
breaks up into a sum of free and trivial $R \langle c \rangle $%
-modules;\ however, by the first part there are no trivial modules in $R [ C_{G} ] $ and the result follows.\ \ \ $\square $

\begin{lemma}\label{lem25}
Again suppose that $c$ cannot be written as a commutator; then
\begin{equation*}
\mathcal{L}\left( 1+ ( 1-c ) I ( R [ G ]  )
\right) \supseteq p\phi \left(  ( 1-c ) I (R[ G ]) \right)
+p^{2}\phi  ( 1-c ) R
\end{equation*}%
and 
\begin{equation*}
\frac{p\phi \left( ( 1-c ) R [ G ] \right) }{\mathcal{L}%
\left( 1+ ( 1-c ) R [ G ] \right) }\cong \frac{R}{\left(
1-F\right) R+pR}.
\end{equation*}
\end{lemma}

Proof. Using the obvious fact that $ ( 1-c ) R [ G ]
=R+I ( R [ G ]  )$, we easily derive the equality 
\begin{equation*}
1+ ( 1-c ) R [ G ] = ( 1+ ( 1-c ) R )\cdot 
 ( 1+p ( 1-c ) R )\cdot   ( 1+ ( 1-c ) I ( R[ G ]  )  ) ;
\end{equation*}%
and so by Lemma \ref{lem20} (a)
\begin{equation*}
\mathcal{L}\left( 1+ ( 1-c ) R [ G ] \right) =\mathcal{L}%
\left( 1+ ( 1-c ) R\ ) +p ( 1-c ) R\ + (
1-c ) \phi  ( I ( R [ G ]  ) \right)
\end{equation*}%
So by Lemma \ref{lem18} (b) we see that 
\begin{equation*}
\frac{p\phi \left( (1-c) R[G] \right) }{\mathcal{L}%
\left( 1+\left( 1-c\right) R[G]  \right) }\cong \frac{\left(
1-c\right) R [ C_{G} ] }{ ( 1-F ) R+pR+ ( 1-c )
\phi( I ( R[G]   ) )}.\ 
\end{equation*}%
By Lemma \ref{lem24}, as $c$ is not a commutator, we have a well-defined map $ (
1-c ) R [ C_{G} ] \rightarrow R/pR$ given by $\phi
 (  ( 1-c ) rg ) \mapsto r{\mod}pR$ for $r\in R$, $g\in
G;$ and the kernel of this map is $ ( 1-c ) pR+ ( 1-c )
\phi \left( I(R[G]) \right) $ and so we have an
isomorphism%
\begin{equation*}
\frac{p\phi \left(  ( 1-c ) R [ G ] \right) }{\mathcal{L}%
\left( 1+ ( 1-c ) R [ G ] \right) } \cong \frac{R}{ (
1-F ) R+pR}.\ \ \ \ \ \square
\end{equation*}

\subsection{Torsion determinants}\label{3.2}

The next result in this section concerns an extension of  Wall's theorem on
torsion determinants for a $p$-adically complete ring $R$. Note that, for
the purposes of this result, we only require that $R$ be a $p$-adically complete
torsion free integral domain.

\begin{thm}\label{thm26}
The group of all elements of finite order in  $\mathrm{Det} ( 1+I ( R[G])  ) $ is equal to the group  $\mathrm{Det} (
G ) $.
\end{thm}

The proof of this result is very similar to that when$\ R$ is a $p$-adic
ring of integers; since the result is important, but the proof is quite
intricate and does not involve any fundamentally new ideas, a detailed proof
is provided as an Appendix in Section 6.\ We can now use this result to show:

\begin{thm}\label{thm27}
The kernel of the map $\nu :\mathrm{Det} ( 1+I ( R[G])  ) \rightarrow p\,\phi ( I ( R[G])    ) $
is equal to  $\mathrm{Det} (
G ) $.
\end{thm}

\begin{thm}\label{thm28}
Let $\nu ^{\prime }$ denote the restriction of  the homomorphism $\nu $ to 
$\mathrm{Det} ( 1+{\mathcal A} ( R[G])  )$; then $\nu^{\prime }$ is an isomorphism 
\begin{equation*}
\mathrm{Det} ( 1+{\mathcal A} ( R[G])  ) \xrightarrow{\ \sim\ }
p\, \phi  ( \mathcal{A}(R[G])) .
\end{equation*}
\end{thm}

{\sl Proof of Theorem \ref{thm27}. }  By Theorem \ref{thm26} it will suffice to show that $\mathrm{Det} ( 1+I ( R[G])  )_{\rm tors} =\ker (\nu)$.
Since ${\rm Im}(\nu)$ is torsion free, it is clear that $\mathrm{Det} ( 1+I ( R[G])  )_{\rm tors} \subset \ker (\nu)$.  We now establish
the opposite inclusion. For this we begin by noting that $ ( 1+I ( R[G])) / ( 1+p^{2}I(R[G])) $ has finite $p$-power exponent, and from Proposition \ref{prop10} we
know that $\upsilon $ is injective on $\mathrm{Det}( 1+p^{2}I ( R[G]))$; thus if $x\in I(R[G])$, then $\mathrm{Det}  ( 1+x ) $ is torsion if and only if 
\begin{equation*}
0=\upsilon  ( \mathrm{Det} ( 1+x )  ) =\phi \circ \log
 ( 1+x ) .
\end{equation*}
Suppose now that $\mathrm{Det}  ( 1+x )  \in \ker (\nu)$; then by the
above it will suffice to show that $\phi \circ \log  ( 1+x ) =0$. 
We prove this by induction on the order of the group $G,\ $noting that, to
start the induction, the result is vacuous if $G$ is trivial. We then adopt
the previous notation, choose a central element of order $p$ and put $%
\overline{G}=G/ \langle c \rangle $ etc. By the induction
hypothesis, $\phi ( \log  ( 1+\overline{x} )  ) =0$ and
so, for any character $\overline{\chi }$ of $G$, which is inflated from $
\overline{G}$, 
\begin{equation*}
0=\overline{\chi }\left( \phi \left( \log \left( 1+x\right) \right) \right).
\end{equation*}%
Recall that the lift of Frobenius $F$ acts on $R[ G] $ and $R[ C_{G}] $ by the coefficients $R$; thus, if $\overline{\chi }$
takes values in the ring $\mathbf{Z}_{p} [ \zeta _{\overline{\chi }}]$, then we have the equality in $R\otimes \mathbf{Z}_{p}[ \zeta
_{\overline{\chi }}]$:
\begin{equation*}
0=F\otimes \mathrm{id}\left( \overline{\chi }\left( \phi \left( \log \left(
1+x\right) \right) \right) \right) =\overline{\chi }\left( \phi \left( \log
\left( 1+Fx\right) \right) \right) .
\end{equation*}%
However, by Lemma 1.8 in Chapter 6 of [T] we know that for any character $%
\chi $ of $G$ the virtual character $\psi ^{p}(\chi) $ is inflated from $%
\overline{G}$. Thus for each character $\chi $ of $G,$ by the second part of
Lemma \ref{lem14}, $\ $we have equalities 
\begin{eqnarray*}
0 &=&\chi \left( \nu \circ \mathrm{Det} ( 1+x ) \right) \\
&=&p\log \left( \mathrm{Det} ( 1+x ) \left( \chi \right) \right)
-\log \left( \mathrm{Det} ( 1+Fx ) \left( \psi ^{p}\chi \right)
\right) \\
&=&p\log \left( \mathrm{Det} ( 1+x ) \left( \chi \right) \right) .
\end{eqnarray*}%
By the first part of Lemma \ref{lem14} this shows that $\phi  ( \log  (
1+x )  ) =0$, as required. \ \ \ $\square \bigskip $

Proof of Theorem \ref{thm28}. By the above
\begin{eqnarray*}
\ker (\nu ^{\prime }) &=&\text{\textrm{Det}} ( 1+\mathcal{A}(R[G])) _{\rm tors} \\
&\subset &\text{\textrm{Det}}( 1+I(R[G]))_{\rm tors}\cap \text{\textrm{Det}} ( 1+\mathcal{A}(R\left[ G\right] ) )
\\
&\subset &\text{\textrm{Det}} ( G ) \cap \text{\textrm{Det}} (
1+\mathcal{A}(R [ G ] ) ) \\
\;\;\;\;\; &\subset &\text{\textrm{Det}} ( \left[ G,G\right]  )
=1.\;\;\;\square
\end{eqnarray*}
\medskip

\subsection{The first exact sequence for determinants.}\label{3.3}

Here we assemble the above to complete the proof of Theorem \ref{thm5} (a).

Using the fact that $\mathrm{Det}$ is an isomorphism on $R[G^{\rm ab}]^\times$, from the above we now deduce that the following
sequence is exact 
\begin{equation}\label{eq16}
0\rightarrow p\,\phi  ( \mathcal{A} ( R [ G ]  ) ) 
\xrightarrow{\nu ^{\prime -1}}\mathrm{Det} (   1+I ( R%
 [ G ]  )   ) \rightarrow 1+I ( R [ G^{\rm ab}]  ) \rightarrow 1.
\end{equation}
 Using the splitting 
 $$
 {\rm Det}
 (  ( 1+R[G] )
 ) = {\rm Det} (   1+I ( R[G])) )
\times R^{\times },
$$ 
we get a further exact sequence
\begin{equation}\label{eq17}
0\rightarrow p\,\phi  ( \mathcal{A} ( R [ G ]  ) ) 
\xrightarrow{\nu ^{\prime -1}} \mathrm{Det} ( R [ G ]
^{\times } ) \rightarrow R[G^{\rm ab}]^\times\rightarrow 1.
\end{equation}
Since $R$ is torsion free, and $\phi(\A(R[G]))$ is a free
$R$-module, $p\, \phi(\A(R[G]))\cong \phi(\A(R[G]))$.
This then completes the proof of Theorem \ref{thm5} (a).
\ \ \ $\square$.
\medskip

\subsection{The second exact sequence; cokernel of the group logarithm}\label{3.4}

Theorem \ref{thm5} (b) now follows from:

\begin{thm}\label{thm29}
Let $G$ be a $p$-group. The following sequence is exact:
\begin{equation}\label{eq18}
\mathrm{Det} ( R [ G ] ^{\times } ) \xrightarrow{\ {p^{-1}\nu }\ }
 \phi  ( I(R [ G ] ) ) \xrightarrow { \ } \frac{R}{
 ( 1-F ) R}\otimes _{\mathbf{Z}_{p}}G^{\rm ab}\rightarrow 0.
\end{equation}
\end{thm}

Proof. We begin by reducing to the case where $G$ is abelian. For this we note
from (\ref{eq17})  and Theorem \ref{thm16} that we have the following commutative
diagram with exact rows:
\begin{equation*}
\begin{array}{ccccccccc}
0 & \rightarrow & \phi  ( \mathcal{A} ) & \rightarrow & \mathrm{Det}%
 ( 1+I(R[G])) & \rightarrow & 1+I(R[G^{\rm ab}])  & \rightarrow & 1 \\ 
&  & \downarrow = &  & \downarrow p^{-1}\nu _{G} &  & \ \ \downarrow p^{-1}\nu
_{G^{\rm ab}} &  &  \\ 
0 & \rightarrow & \phi ( \mathcal{A} ) & \rightarrow & \phi  (
I(R[G])) & \rightarrow & I(R[G^{\rm ab}]) & 
\rightarrow & 0.%
\end{array}%
\end{equation*}%
So by the Snake lemma we conclude that $\mathrm{coker} ( p^{-1}\nu
_{G} ) \cong \mathrm{coker} ( p^{-1}\nu _{G^{\rm ab}} ) $, which
indeed means that we may now restrict our attention to the case when $G$ is
abelian.

By Proposition \ref{prop23} we know that $(p^{-1}\nu _{G}) ( 1+ I^2(R[G ])) =I^2(R[G ])$ and so to calculate the cokernel of 
$p^{-1}\nu _{G}$, it will suffice to determine $(p^{-1}\nu _{G}) ( 1+ I(R[G ])) \mod I^2(R[G ])$.
\medskip

Next we briefly recall the concept of the cotangent space at the origin for
a commutative augmented $R$-algebra $A$. For such an algebra we write $\varepsilon $ for the augmentation map and $I_{A}$ for the augmentation
ideal; so that we have a direct decomposition of $R$-modules $A=R\oplus
I_{A}$ and the cotangent space at the origin is%
\begin{equation*}
T_{A}^{\ast }=I_{A}/I_{A}^{2}.
\end{equation*}%
For an element $f\in I_{A}$ we shall write $df$ for its image in $
T_{A}^{\ast }$. For future reference note that 
\begin{equation}\label{eq19}
d\log  ( 1+f ) =\frac{df}{1+f}=df.
\end{equation}

Suppose now that we can write the abelian $p$-group $G$ as 
\begin{equation*}
G=\bigoplus _{i=1}^{n} {\mathbf{Z}}/{p^{m_{i}}\mathbf{Z}}
\end{equation*}%
and we let $g_{i} \in G$ be chosen with zero component away from $i$ and
to have $i$-th component a generator $\mathbf{Z}/p^{m_{i}}\Z$.
Let $R\left[ T\right] $ denote the polynomial $ R$-algebra $R [ T_{1},\ldots
,T_{n} ] $ of (commuting) indeterminates $T_{i}$. We write $\theta :R[ T ] \rightarrow R [ G ] $ for the $R$-algebra
homomorphism with the property that $\theta ( T_{i} ) =g_{i}$. We
then have the following presentation of the group ring $R\left[ G\right] :$%
\begin{equation*}
0\rightarrow N\rightarrow R [ T ] \overset{\theta }{\rightarrow }R%
 [ G ] \rightarrow 0
\end{equation*}%
where $N$ denotes the $R [ T ] $-ideal generated by the terms $%
T_{i}^{p^{m_{i}}}-1$ for $i=1,\ldots ,n$; and we have the following
equalities: 
\begin{eqnarray*}
T_{R\left[ T\right] }^{\ast } &=&\oplus _{i=1}^{n}R\cdot d ( T_{i}-1 ) ,
\\
dN &=&\oplus _{i=1}^{n}p^{m_{i}}R\cdot d ( T_{i}-1 ),
 \\
T_{R [ G ] }^{\ast } &=&\oplus _{i=1}^{n}\frac{R}{p^{m_{i}}R}\cdot 
d ( g_{i}-1 ) .
\end{eqnarray*}

Recall that $R$ supports a $\mathbf{Z}_{p}$-algebra homomorphism  $F$ (the
lift of Frobenius); we now vary from the convention established in the
Introduction and extend $ F$ to $R [ T ] $ and to $R [ G ]
$ by setting $F ( T_{i} ) =T_{i}^{p}$ and $F ( g_{i} )
=g_{i}^{p}$ for all $i$. The following result should properly be seen in the
context of the action of Frobenius on differentials. An \textit{ad hoc}
treatment is provided here; for a full treatment the reader is referred to
[I].

\begin{prop}\label{prop30}
(a) There is an $F$-semi-linear homomorphism, also denoted $F$ of $T_{R [ T%
 ] }^{\ast }$ with the property that $pF\circ d=d\circ F$ and with $%
F ( r\cdot d(T_{i}-1) ) =F ( r )\cdot  d(T_{i}-1)$, for $r\in R $; that is to
say, since $R$ has characteristic zero and since $T_{R [ T ]
}^{\ast }$ is $R$-free, for $f\in I_{R [ T ] }$, we have the
equality 
\begin{equation}\label{eq20}
F ( df ) =\frac{1}{p}\cdot d ( F ( f ) ) .
\end{equation}

(b) $F$ maps $dN$ into itself, and so induces an endomorphism of $T_{R[ G ] }^{\ast }$. Recall that for $j\in I_{R [ G ] }$,
\begin{equation*}
p^{-1}\mathcal{L}\left( 1+j\right) =\log \left( 1+j\right) -p^{-1}\log
\left( 1+F\left( j\right) \right) \in I_{R [ G ] }.
\end{equation*}

(c) For $j\in I_{R [ G ] }$ we have   
$$
d\left( p^{-1}\mathcal{L}\left( 1+j\right) \right) =\left( 1-F\right) dj
$$
and so 
$$
d\left( p^{-1}\mathcal{L}\left( 1+I_{R [ G ] }\right)
\right) =\left( 1-F\right) T_{R [ G ] }^{\ast }.
$$
\end{prop}

Proof. (a)\ For $r\in R$, we set $F ( r\cdot d ( T_{i}-1 )  )
=F ( r )\cdot d ( T_{i}-1 )$; we then need only note that 
\begin{eqnarray*}
d ( F ( r\cdot ( T_{i}-1 )  )  ) &=&d ( F (
r )\cdot  ( T_{i}^{p}-1 )  ) =F ( r )\cdot  d (  (
T_{i}^{p}-1 )  ) \\
&=&p\, F ( r)\cdot T_{i}^{p-1}d (  ( T_{i}-1 )  )
=p F ( r )\cdot d (  ( T_{i}-1 )  ) .
\end{eqnarray*}%

To prove (b) we choose $f\in R [ T ] $ with $\theta (
f ) =j$. Then from (a) and (\ref{eq18}) above we know that
\begin{equation*}
d ( \log \left( 1+j\right) -p^{-1}\log \left( 1+F ( j )  )) =dj-p^{-1}d ( F( j ) \right) = ( 1-F ) dj.
\end{equation*}%
Applying $\theta $ then gives the desired equality. Part (c) now follows.\ \ \ $\square$
\bigskip

We now complete the proof of Theorem \ref{thm29}. Recall that we have reduced to the
case in which $ G=G^{\rm ab}$ is abelian. We write $T_{R [ G ]
}^{\ast }=R\otimes T_{\mathbf{Z}_{p} [ G ] }^{\ast }$, and note
that, since we have $(1-F) ( R\otimes T_{\mathbf{Z}_{p} [ G ]}^{\ast
}) = ( 1-F )R\otimes T_{\mathbf{Z}_{p}\left[ G\right]
}^{\ast },$ it now follows that 
\begin{equation*}
d ( p^{-1}\mathcal{L} ( 1+I_{R [ G ] } )  )
= ( 1-F ) T_{R\left[ G\right] }^{\ast }= ( 1-F ) R\otimes
T_{\mathbf{Z}_{p} [ G ] }^{\ast }
\end{equation*}%
and 
\begin{equation*}
\frac{T_{R [ G ] }^{\ast }}{ ( 1-F ) T_{R [ G ]
}^{\ast }}=\frac{R\otimes T_{\mathbf{Z}_{p} [ G ] }^{\ast }}{ (
1-F ) R\otimes T_{\mathbf{Z}_{p} [ G ] }^{\ast }}=\frac{R}{%
 ( 1-F ) R}\otimes G
\end{equation*}%
as required. \ \ $\square$ 
\medskip 

As an immediate application we now prove:

\begin{lemma}
Let $c$ be a central element of order $p$ with $c\in  [ G,G ]$,
and, as previously, we put $\overline{G}=G/ \langle c \rangle $.
Then the following map, induced by the group logarithm $\nu$, is an
isomorphism:
\begin{equation*}
\ker  ( \mathrm{Det} ( 1+I ( R [ G ]  )  )
\rightarrow \mathrm{Det} ( 1+I ( R [ \overline{G} ]  )
 )  ) \rightarrow \ker  ( \phi  ( I ( R [ G ]
 )  ) \rightarrow \overline{\phi } ( I ( R [ \overline{
G} ]  )  )  ) .
\end{equation*}
\end{lemma}

Proof. We set ${\mathcal K}_{1}=\ker  \left[ \mathrm{Det} ( 1+I ( R [ G
 ] )  ) \rightarrow \mathrm{Det} ( 1+I ( R [ 
\overline{G} ]  )  )  \right] $ and similarly, ${\mathcal K}_{2}=\ker\left[
\phi  ( I ( R [ G ]  )  ) \rightarrow \overline{
\phi } ( I ( R [ \overline{G} ] ) )  \right] $
 and consider the diagram 
\begin{equation*}
\begin{array}{ccccccccc}
1 & \rightarrow & {\mathcal K}_{1} & \rightarrow & \mathrm{Det}\left( 1+I ( R [
G ]  ) \right) & \rightarrow & \mathrm{Det}\left( 1+I ( R [
\overline{G} ]  ) \right) & \rightarrow & 1 \\ 
&  & \downarrow &  & \downarrow &  & \downarrow &  &  \\ 
0 & \rightarrow & {\mathcal K}_{2} & \rightarrow & \phi \left( I ( R [ G ]
 ) \right) & \rightarrow & \phi \left( I ( R [ \overline{G}
 ]  ) \right) & \rightarrow & 0%
\end{array}%
\end{equation*}%
where the vertical maps are induced by the group logarithm $\nu$. Since $G^{\rm ab}=\overline{G}^{\rm ab}$, we know that the central and right-hand vertical
maps have the same kernel (by Theorem \ref{thm26}) and the same cokernel (by Theorem
\ref{thm29}); therefore, by the snake lemma, we deduce that the left-hand vertical
map is an isomorphism.\ \ \ \ $\square $
\bigskip

\section{Descent and fixed points for $p$-groups}\label{descent}
\setcounter{equation}{0}

\subsection{Fixed point theorem for $p$-groups.}\label{3.5}

In this subsection we adopt the notation of Theorem \ref{thm6}. Thus $S$ is an
extension ring of $R$, which is a finitely generated, free $R$-module; so
that $S$ is also $p$-adically complete. We suppose that $S$ satisfies the
remaining conditions in the Hypothesis stated in the Introduction; so that $S$ is an integral domain of characteristic zero, which supports a lift of
Frobenius. Note that we do not need to suppose that  the lift of
Frobenius $F_{R}$ is compatible with the lift of Frobenius $F_{S}$, so that 
${F_{S}} {|R}$ need not equal $F_{R}$. Note also that at this point we do
not need to suppose that these two rings are Noetherian and normal. We
suppose that the field of fractions of $S$ is Galois over $N$ with Galois
group $\Delta $ acting on $S$ with $S^{\Delta }=R$. Then, since $S$ is
finite over $R$, we know that $(S^{\times})^{\Delta }=R^{\times }$.

The main result of this section is that, with the above notation and
hypotheses,

\begin{thm}\label{thm32}
We have the equality
\begin{equation*}
\mathrm{Det}\left( 1+I ( S [ G ]  )  \right) ^{\Delta }=%
\mathrm{Det}\left( 1+I ( R [ G ]  )  \right) ;
\end{equation*}%
hence using the commutative diagram 
\begin{equation*}
\begin{array}{ccccccc}
1\rightarrow & \mathrm{Det}\left( 1+I ( S [ G ]  ) \right)
^{\Delta } & \rightarrow & \mathrm{Det} ( S [ G ] ^{\times
} ) ^{\Delta } & \rightarrow & (S^{\times})^{ \Delta } & \rightarrow 1 \\ 
& \uparrow &  & \uparrow &  & \uparrow = &  \\ 
1\rightarrow & \mathrm{Det}\left( 1+I ( R [ G ]  ) \right)
& \rightarrow & \mathrm{Det} ( R [ G ] ^{\times } ) & 
\rightarrow & R^{\times } & \rightarrow 1%
\end{array}%
\end{equation*}%
we see that ${\rm Det} ( S [ G ] ^{\times } ) ^{\Delta
}={\rm Det} ( R [ G\ ] ^{\times } )$.
\end{thm}
\medskip

Proof. We argue by induction on the order of $G$, and we note that the
result is clear if $G$ is abelian, since in that case we have 
\begin{equation*}
\mathrm{Det}( S[ G] ^{\times }) ^{\Delta }\cong (S[
G] ^{\times})^{ \Delta }=R[ G] ^{\times }\cong \mathrm{Det}
( R[ G] ^{\times }) .
\end{equation*}
We therefore suppose that $G$ is non-abelian. As previously we select a
central commutator $c\in G$ of order $p$, put $\overline{G}=G/\langle
c\rangle$, and write $q$ for the quotient map $q:G\rightarrow 
\overline{G}$; then, by the induction hypothesis, the theorem is true for $%
\overline{G}$. We then have the commutative diagram with exact rows 
\begin{equation*}
\begin{array}{ccccccc}
0\rightarrow & p\phi ( \mathcal{A}( R[ G] )
 ) & \rightarrow & \mathrm{Det}( 1+I( R[ G]
)) & \rightarrow & 1+I (R[ G^{\rm ab}] ) & 
\rightarrow 1 \\ 
& \downarrow &  & \downarrow q &  & \downarrow = &  \\ 
0\rightarrow & p\overline{\phi } ( \mathcal{A} ( R [ \overline{G}]))& \rightarrow & \mathrm{Det}( 1+I( R[ 
\overline{G}]))& \rightarrow & 1+I(R[
\overline{G}^{\rm ab}])& \rightarrow 1.
\end{array}
\end{equation*}

\begin{lemma}\label{lem33}
The map $\nu$ of Corollary \ref{cor15} induces an isomorphism
\begin{equation}\label{eq21}
\ker\left[\mathrm{Det} ( 1+I( R[ G] ))
\rightarrow \mathrm{Det}(1+I( R[ \overline{G}])\right]\cong \mathrm{Det}( 1+( 1-c) R[ G]) .
\end{equation}
\end{lemma}

Proof. By the snake lemma we have the equality
\begin{equation}\label{eq22}
\ker \left( \mathcal{A} ( R [ G ] ) \rightarrow \mathcal{A%
} ( R [ \overline{G} ] ) \right) =\ker\left(R[ G] \rightarrow R [ \overline{G} ] \right) = ( 1-c ) R[ G ] .
\end{equation}%
Applying the snake lemma again we have%
\begin{equation*}
\ker \left[\mathrm{Det}(1+I (R[G]))
\rightarrow \mathrm{Det}(1+I(R[\overline G]))\right]=\ker \left[\mathrm{Det} ( 1+\mathcal{A}(R[G])) \rightarrow \mathrm{Det}(1+\mathcal{A}(R[\overline G])\right]
\end{equation*}
and by Theorem \ref{thm28} this is isomorphic to 
\begin{equation*}
\ker \left[ p\,\phi ( \mathcal{A} (R[G])
\rightarrow p\,\phi  ( \mathcal{A} (R[\overline{G}]) \right] =p\, \phi  (  ( 1-c ) R [ G ]
 ) \cong \mathrm{Det} ( 1+ ( 1-c ) R [ G ] )
\end{equation*}
  with the latter isomorphism coming from Lemma \ref{lem21}.\ \ \ $\square $
  \medskip

We now complete the proof of the Theorem \ref{thm32}. Suppose we are given 
$${\rm Det}(x) \in \mathrm{Det} ( 1+I ( S [ G ]))^{\Delta }$$ and we put $\overline{x}=q ( x )$.
Then since the theorem is true for $\overline{G}$
\begin{equation*}
{\rm Det} ( \overline{x} ) \in \mathrm{Det} (
1+I ( S [ \overline{G} ]  )  ) ^{\Delta }=\mathrm{Det}%
 ( 1+I ( R [ \overline{G} ] )  ) \text{ }
\end{equation*}
and, because $I ( R [ G ]  ) $ maps onto $I ( R [ 
\overline{G} ]  )$, we know that $q$ maps $\mathrm{Det} (
1+I ( R [ G ]  )  ) $ onto ${\rm Det} (
1+I ( R [ \overline{G} ]  )  )$. We can therefore
find $y\in 1+I ( R [ G ]  ) $ such that $q\left( \mathrm{%
Det}\left( xy^{-1}\right) \right) =1.$ By Lemma \ref{lem33}, with $ R$ replaced by $S$, we know that $\mathrm{Det} ( xy^{-1} ) \in \mathrm{Det} (
1+ ( 1-c ) S [ G ]  ) ^{\Delta }$. Thus it will now
suffice to show that 
\begin{equation}\label{eq23}
\mathrm{Det} ( 1+ ( 1-c ) S [ G ]  ) ^{\Delta }=%
\mathrm{Det} ( 1+ ( 1-c ) R [ G\ ]  ) .
\end{equation}
Unfortunately we cannot deduce this immediately, since the map $\nu $ does
not, in general, respect the action of $\Delta $ (due to the appearance of
the $F$-semilinear map $\Psi $).  However, we know from Lemma \ref{lem19} that the
map $\Psi $ is null on $ ( 1-c ) S [ G ] $ and so the
restriction of  $\mathcal{L}$ to $1+ ( 1-c ) S [ G ] $
is just $p$ times the usual $p$-adic logarithm of Proposition \ref{prop10}; thus we
see that 
\begin{equation*}
{\mathcal L}_{S}\mid_{\,1+( 1-c) R[ G] }=\mathcal{L}_{R}\mid _{\, 1+( 1-c) R[ G] }
\end{equation*}
 and so 
\begin{equation*}
\nu_{S}\mid_{\,\mathrm{Det}( 1+( 1-c) R[ G]) }=\nu _{R}\mid _{\, \mathrm{Det}( 1+( 1-c)R[ G]) }
\end{equation*}%
and that furthermore $\mathcal{L}_{S}$, resp. $\nu _{S}$, both commute with $\Delta $-action \textit{when restricted to} $1+ ( 1-c ) S [ G%
 ] $, resp. ${\rm Det}( 1+ ( 1-c ) S [ G ])$. Thus by Lemma \ref{lem22} we get 
\begin{eqnarray*}
\nu_{S}\left( \mathrm{Det} ( 1+ ( 1-c ) S [ G ]
 )^{\Delta } \right) &\subset & \left( p\phi  (  ( 1-c ) S
 [ G ]  )  \right) ^{\Delta } \\
&=&p\,\phi  (  ( 1-c ) R [ G ]  ) \\
&=&\nu _{R} \left( \mathrm{Det} ( 1+ ( 1-c ) R [ G ]
 )  \right) \\
&=&\nu _{S} \left( \mathrm{Det} ( 1+ ( 1-c ) R [ G\ ]
 )  \right) \;
\end{eqnarray*}%
and the result follows since, by Theorem \ref{thm28}, $\nu _{S}$ is injective on 
${\rm Det}( 1+\mathcal{A} ( S [ G ]  )  )$
. \ \ \ \ \ $\square$ 
\medskip  

\subsection{Norm maps for $p$-groups.}\label{3.6}

The notation and hypotheses are as in the previous subsection; in
particular, we again impose the hypotheses of Theorem \ref{thm6}. Thus $S$ is a ring
extension of $R$, which is a finitely generated, free $R$-module,
and we again suppose that both $S$ and $R$ satisfy the standing hypothesis.
We also suppose that the field of fractions of $S$ is a Galois extension of $N$ with Galois group $\Delta$, which acts on $S$ with the property that $S^{\Delta }=R$; in addition we now also suppose that  the trace map ${\rm Tr}:S\rightarrow R$ is surjective. Recall that, by taking bases, we have the
co-restriction (or Fr\"{o}hlich norm) map (see for instance 1.4 in [T]) 
\begin{equation*}
\mathcal{N}:\mathrm{Det}( S[ G] ^{\times })
\rightarrow \mathrm{Det}( R[ G] ^{\times })
\end{equation*}
which is characterized by the property%
\begin{equation}
\mathcal{N}( \mathrm{Det}( x) ) =\prod\limits_{\delta
\in \Delta }\mathrm{Det}( x^{\delta }) =\mathrm{Det}\left(
\prod\limits_{\delta \in \Delta }x^{\delta }\right) .
\end{equation}
By considering the augmentation map $G\rightarrow \left\{ 1\right\} $ and
the abelianisation map $G\rightarrow G^{\rm ab},$ it follows at once from the
definitions that we have the inclusions
\begin{eqnarray*}
\mathcal{N} \left( \mathrm{Det} ( 1+I(S[G]))  \right) &\subset &\mathrm{Det} ( 1+I ( R [ G ]
 )  ) \\
\mathcal{N}\left(\mathrm{Det} ( 1+\mathcal{A} ( S [ G ]  )
 ) \right) &\subset &\mathrm{Det} ( 1+\mathcal{A} ( R [ G ]
 )  ) .
\end{eqnarray*}

With the above notation and hypotheses

\begin{thm}\label{thm34}
We have equalities 
\begin{equation}\label{eq25}
\mathcal{N}\left( \mathrm{Det} ( 1+I ( S [ G ]  )
 ) \right) =\mathrm{Det} ( 1+I ( R [ G ]  )
 )
\end{equation}
\begin{equation}\label{eq26}
\mathcal{N}\left( \mathrm{Det} ( 1+\mathcal{A} ( S [ G ]
 )  ) \right) =\mathrm{Det} ( 1+\mathcal{A} ( R [ G
 ]  )  ) ;
\end{equation}
while  the full norm map $\mathcal{N}:\mathrm{Det} ( S [ G ]
^{\times } ) \rightarrow \mathrm{Det} ( R[ G ] ^{\times
} ) $ is surjective if and only if we have $N_{S/R} ( S^{\times } )
=R^{\times }$.
\end{thm}

Proof. We first prove (\ref{eq25}) and (\ref{eq26}); we then conclude by considering the
full norm map.

{\sl Step 1.}\ We first suppose that $G$ is abelian. For brevity we write $%
I_{S}=I ( S [ G ]  )$, $I_{R}=I ( R [ G ]
 ) $ and we filter $1+I_{S}$, resp. $1+I_{R}$, by the subgroups $\left\{ \left(
1+I_{S}^{n}\right) \right\} _{n\geq 1}$, resp. $\left\{ \left(
1+I_{R}^{n}\right) \right\} _{n\geq 1}$. By use of successive
approximation, the surjectivity of the trace, and the commutative squares: 
\begin{equation*}
\begin{array}{ccc}
(1+I_{S}^{n}) / (1+I_{S}^{n+1}) & \rightarrow & I_{S}^{n}/
I_{S}^{n+1}=\left( I_{R}^{n}/I_{R}^{n+1}\right) \otimes _{R}S \\ 
\mathcal{N}\downarrow &  & \downarrow {\rm Tr} \\ 
(1+I_{R}^{n})/(1+I_{R}^{n+1}) & \rightarrow & \left( I_{R}^{n}/
I_{R}^{n+1}\right) .%
\end{array}%
\;
\end{equation*}%
together with a standard convergence argument, it then follows that $
\mathcal{N}$ maps $1+I_{S}$ onto $1+I_{R}$.
\smallskip

{\sl Step 2.}\  Consider $\mathrm{Det}(x) \in \mathrm{Det} (
1+I ( R [ G ]  )  ) $ and let $\widetilde{x}$ denote
the image of $x$ in $R [ G^{\rm ab} ]$. By the first step, we can
find $\widetilde{y}\in 1+I ( S [ G^{\rm ab} ]  ) $ such that $%
\mathcal{N} ( \mathrm{Det} ( \widetilde{y} )  ) =\mathrm{%
Det} ( \widetilde{x} ) $ and we can then choose $y\in 1+I ( S%
 [ G ]  )  $ with image equal to $\widetilde{y}$ in $ S [
G^{\rm ab} ]$.  Then $\mathrm{Det}( x) \mathcal{N}( 
\mathrm{Det}( \widetilde{y})) ^{-1}\in \mathrm{Det}(
1+\mathcal{A}( R[ G])) $ and so we are now
reduced to showing 
\begin{equation}\label{eq27}
\mathcal{N}\left( \mathrm{Det}( 1+\mathcal{A} ( S [ G ]
 )  ) \right) \supseteq \mathrm{Det} ( 1+\mathcal{A} ( R[ G] )) .
\end{equation}%
To prove the inclusion (\ref{eq27}) we argue by induction on the order of the
commutator group $[ G,G]$. The inclusion is vacuous when $G$ is
abelian; so this starts our induction, and we now adopt the notation of the
previous subsection; in particular we choose a central commutator $c$ of
order $p$, we put $\overline{G}=G/\langle c\rangle $ and write $q$
for the natural quotient map from $G$ to $\overline{G}$. Let $\mathrm{Det}(x) \in \mathrm{Det}( 1+\mathcal{A}( R[ G]))$. Then, by the induction hypothesis, we can find 
$y\in 1+I(S[G])$ such that $\mathrm{Det}( q(
x)) =\mathcal{N}(\mathrm{Det}(q(y))$; hence by (\ref{eq22}) in the proof of Lemma 
\ref{lem33} 
\begin{eqnarray*}
\mathrm{Det} ( x ) \mathcal{N} ( \mathrm{Det} ( y ))^{-1} &\in &\ker \left[ \mathrm{Det}( 1+\mathcal{A}( R[ G])) \rightarrow \mathrm{Det}( 1+\mathcal{A}( R[ \overline{G}] )) \right] \\
&=&\mathrm{Det}(1+(1-c)R[G]).
\end{eqnarray*}%
It will therefore suffice to show: 
\begin{equation}
\mathcal{N}\left( \mathrm{Det}( 1+( 1-c) S[ G]
) \right) =\mathrm{Det}( 1+( 1-c) R[ G]).
\end{equation}%
As in the previous section we know that $\Psi $ is null on $(
1-c) S[G]$; hence the restriction of $\nu_{S}$ to $\mathrm{Det}( 1+( 1-c) S[ G]) $ is a map 
of $\Delta$-modules and furthermore 
\begin{equation*}
\nu _{S}\mid _{\,\mathrm{Det}( 1+( 1-c) R[ G]) }=\nu _{R}\mid_{\,\mathrm{Det}( 1+( 1-c) R[ G]) },
\end{equation*}%
as both sides agree with $p\,\phi \circ \log$. Hence the following diagram
commutes:
\begin{equation*}
\begin{array}{ccccc}
\mathrm{Det}\left( 1+( 1-c) S[ G] \right) & \xrightarrow{\overset{\nu_S}{\sim}} & p\,\phi  (  ( 1-c ) S [ G ]  )
& = & p\,\phi  (  ( 1-c ) R [ G ]  ) \otimes _{R}S
\\ 
\downarrow \mathcal{N} &  & \downarrow {\rm Tr} &  &  \\ 
\mathrm{Det}\left( 1+ ( 1-c ) R [ G ] \right) & \xrightarrow{\overset{\nu_R}{\sim}}  & p\,\phi  (  ( 1-c ) R [ G ]  )
&  & 
\end{array}%
\end{equation*}%
and the surjectivity of the left-hand vertical arrow follows from the
surjectivity of the trace from $S$ onto $R$.

To prove the final part of the theorem, we use the following commutative
diagram with exact rows: 
\begin{equation*}
\begin{array}{ccccccc}
1\rightarrow & \mathrm{Det}\left( 1+I(S[G]) \right)
& \rightarrow & \mathrm{Det} ( S [ G ] ^{\times } ) & 
\rightarrow & S^{\times } & \rightarrow 1 \\ 
& \downarrow \mathcal{N} &  & \downarrow \mathcal{N} &  & \ \ \downarrow N &  \\ 
1\rightarrow & \mathrm{Det}\left( 1+I ( R [ G ]  ) \right)
& \rightarrow & \mathrm{Det} ( R [ G ] ^{\times } ) & 
\rightarrow & R^{\times } & \ \rightarrow 1.\;%
\end{array}%
\end{equation*}
$\square$
\bigskip

\section{Character action and reduction to elementary groups}\label{4}
\setcounter{equation}{0}

\subsection{Character action on $\Kr_{1}$}\label{4.1}

In this section we do not require that $R$ supports a lift of Frobenius.
Let $\Gr_{0}(\Z_p[G]) $ denote the
Grothendieck group of finitely generated $\Z_p[G]$-modules and let
$\Gr^{\Z_p}_{0}(\Z_p[G]) $  denote the Grothendieck group of finitely generated $\Z_p[G]$-modules  which are projective over $\Z_{p}$. From
38.42 and 39.10 in [CR2] we have:

\begin{prop}\label{prop35}
There are isomorphisms%
\begin{equation*}
\Gr^{\Z_p}_{0}(\Z_p[G]) \cong \Gr_{0}(\Z_p[G])\cong \Gr_{0}(\Q_p[G])
\end{equation*}%
with the first isomorphism induced by inclusion of categories and the second
isomorphism induced by the extension of scalars map $\otimes_{\Z_p}\Q_p$.
\end{prop}

\begin{prop}\label{prop36}
The ring $\Gr^{\Z_p}_{0}(\Z_p[G])$ and hence, by the previous proposition, $\Gr_{0}(\Q_p[G])$, acts naturally on $\Kr_{1}(R[G]) $
via the rule that for an $\Z_p[G]$-lattice $L$ and
for an element of $\Kr_{1}(R[G]) $ represented by a
pair $ ( P,\alpha  ) $ (where $P$ is a projective $R [ G ]$-module and $\alpha $ is an $R [ G ] $-automorphism of $P$), then 
\begin{equation*}
L \cdot  ( P,\alpha  ) = (  ( L\otimes _{\Z_p}P )
,\,  ( 1\otimes _{\Z_p}\alpha  ) ) .
\end{equation*}%
The functor $G\mapsto  \Kr_{1}(R[G]) $ is a
Frobenius module for the Frobenius functor $G\mapsto \Gr^{\Z_p}_{0}(\Z_p[G])$
(see page 4 in [CR2]
and also [L]). The group $\SK_{1} ( R [ G ]  ) $ is a
Frobenius submodule of $\Kr_{1} ( R [ G ]  ) $. Therefore  
the action of $\Gr_{0}(\Q_p[G])$ on $\Kr_{1}(R[G])$ also induces an action on $\mathrm{Det}%
 ( \GL(R [ G ] ) ) $ (see Ullom's Theorem in 2.1 of [T],
and see also below for his explicit description of this action)).
\end{prop}

Proof. From Ex 39.5 in [CR2] we know that $G\mapsto  \Kr_{1}(\Z_p[G]) $ is a Frobenius module for $G\mapsto \Gr^{\Z_p}_{0}(\Z_p[G])$. 
Moreover the extension of scalars map $\Gr^{\Z_p}_{0}(\Z_p[G])\rightarrow 
\Gr^{R}_{0}(R[G])$
is a morphism of Frobenius modules over $\Gr^{\Z_p}_{0}(\Z_p[G])$   by 38.11 loc. cit..
The first part of the proposition follows. Because $\Gr^{R}_{0}(R[G])
\rightarrow  \Gr^{N^c}_{0}(N^c[G])$
  is similarly a morphism of Frobenius modules
over $\Gr^{R}_{0}(R[G])$, it follows that the kernel 
\begin{equation*}
\SK_{1} ( R[G]) =\ker \left[ \Kr_{1}(R[G])  \rightarrow \Kr_{1}(N^c[G]) \right]
\end{equation*}%
is also a Frobenius module over $\Gr^{\Z_p}_{0}(\Z_p[G])$.\ \ \ \ $\square $
\medskip 

Next we recall Ullom's explicit formula for the action of the character ring 
$\Gr_{0}(\Q_p[G])$ on $\mathrm{Det} (\GL ( R [ G ] )  )$.

We view $\Gr_{0}(\Q_p[G])$  as the ring of
characters of finitely generated $\Q_{p}[G]$-modules
and we let $\phi \in \Kr_{0}( \Q_{p}^{c}[ G] )$
(which we identify with the ring of virtual  $\Q_{p}^{c}$-valued
characters of $G$). We consider a pair $ ( P,\alpha  ) $, as above,
and suppose that $P\oplus Q=\oplus _{i=1}^{n}R[ G]$. Then $ ( P,\alpha  ) $ is the element of $\Kr_{1}(R[G])$ given by the class $r=\alpha \oplus {\rm id}_{Q}\in \GL_{n}(R[G])$ 
(see for instance 3.1.7 in [R]). Then Ullom has shown
that the induced Frobenius action of $\Gr_{0}(\Q_p[G])$ on $\mathrm{Det} ( \Kr_{1} (R[G]))$ is given explicitly by 
\begin{equation*}
\theta \cdot ( P,\alpha  ) =\theta \cdot  ( \phi \mapsto \mathrm{Det}%
 ( r )  ( \phi  ) ) = ( \phi \mapsto \mathrm{Det}(r)  ( \overline{\theta }\phi  )  )
\end{equation*}
where $\overline{\theta }$ denotes the contragredient of the character $\theta$. In particular we have the identity 
\begin{equation}\label{eq29}
 ( \mathrm{Ind}_{H}^{G}\theta  ) \cdot \mathrm{Det} ( r ) =%
\text{\textrm{Ind}}_{H}^{G}\left( \theta \cdot  ( \text{\textrm{Res}}%
_{G}^{H} ( \mathrm{Det}(r)  )  ) \right)
\end{equation}%
which is one of the standard identities for Frobenius modules. The proofs of
these standard facts are exactly the same as the proofs in Chapter 2 pages
21-25 in [T].

\begin{prop}\label{prop37}
For each $n\geq 1$, $\Kr_{1} ( R [ G ], p^{n} ) $ is a $\Gr^{\Z_p}_{0}(\Z_p[G])$-Frobenius submodule of $\Kr_{1}(\Z_p [ G ] )$,
and so the action of $\Gr^{\Z_p}_{0}(\Z_p[G])$ on 
$\SK_{1}(R[G])$ induces a
Frobenius module structure on $\widehat{\SK}_{1}(R[G])$.
\end{prop}

Proof. The first part follows from the exact sequence%
\begin{equation*}
0\rightarrow \Kr_{1} ( R [ G ], p^{n} )  \rightarrow
\Kr_{1} ( R [ G ] )  \rightarrow \Kr_{1} ( R_n [ G ] )  \rightarrow 0
\end{equation*}%
and the fact that the extension of scalars map $\Gr^{R}_{0}(R[G])\rightarrow 
\Gr^{R_n}_{0}(R_n[G])$ is a morphism of Frobenius modules, so that the kernel in the above exact
sequence is also a $\Gr^{R}_{0}(R[G])$ Frobenius
submodule of $\Kr_{1} ( R [ G ] )$. The second part then
follows since the image of $\SK_{1} ( R [ G ] )$ in $\Kr_{1} ( R_n [ G ] )$
  is a Frobenius module, because it
is naturally isomorphic to 
\begin{equation*}
\frac{\SK_{1}(R[G]) +\Kr_{1}( R[ G]
, p^{n}) }{\Kr_{1}( R[ G] , p^{n}) }.\ \ \ \ \square
\end{equation*}
\smallskip

\subsection{Brauer Induction}

Let $\mu _{m}$ denote the group of roots of  unity of order $m$ in $\Q_p^{c}$. We then identify $\Gal (\Q_p(\mu_m)/\Q_p) $ as a subgroup of $ (\Z/m\Z)^{\times }$ in the usual way. Recall that a
semi-direct product of a cyclic group $C$ (of order $m$, say, which is
coprime to $p$) by a $p$-group $P$, $C\rtimes P,$ is called $\Q_{p}$-$p$-elementary (see page 112 in [S]) if, for each $\pi \in P$, there exists 
\begin{equation*}
t=t( \pi ) \in \Gal (\Q_p(\mu_m)/\Q_p)  \subseteq  (\Z/m\Z)^{\times }
\end{equation*}%
such that for all $c\in C$ 
\begin{equation*}
\pi c\pi ^{-1}=c^{t}.
\end{equation*}

\begin{thm}\label{thm38}
For a given finite group $G$, there exists an integer $l$ coprime to $p$
such that 
\begin{equation*}
l\Gr_{0} (\Q_p[G]) \subseteq \sum\limits_{J}
\mathrm{Ind}_{J}^{G}\left( \Gr_{0} ( \Q_{p} [ J ]  )
\right)
\end{equation*}
where $J$ ranges over the $\Q_{p}$-$p$-elementary subgroups of $G$.
\end{thm}

Proof. See Theorem 28 in [S].\ \ \ \ \ $\square $
\bigskip

\section{Determinants for elementary groups.}\label{5}
\setcounter{equation}{0}

Throughout this section we shall suppose that, in addition to the standing
hypotheses, $R$ is also a normal Noetherian ring. We begin with an
algebraic result which we shall require later in this section.

Suppose that $\O_{K}$ is the ring of integers of a finite non-ramified
extension $K$ of $\Q_{p}$. Recall that $N$ is the field of
fractions of $R$. Let $S$ denote the ring $R\otimes_{\Z_{p}}\O_{K}$
and let $M$ be the ring of fractions of $S$. Since $M$ is a separable $N$-algebra, 
it can be written as a finite product of field extensions $M_{i}$ of $N$:
\begin{equation}\label{eq30}
M=\prod\limits_{i=1}^{n}M_{i}.
\end{equation}%
Since $K$ is a finite non-ramified extension of $\Q_{p}$, we know
that $S$ is etale over $R$ and hence is normal (see e.g. page 27 in [Mi]).
If $S_{i}$ is the normalization of $R$ in $M_{i}$, then
\begin{equation}\label{eq31}
S=\prod\limits_{i=1}^{n}S_{i}.
\end{equation}%

\begin{lemma}\label{lem39}
Let $F$ denote the lift of Frobenius on $S$ given by the tensor product of
the lift of Frobenius of $R$ with the Frobenius automorphism of $\O_{K}$;
then $F( S_{i}) \subset S_{i}$.
\end{lemma}

Proof. Let $ \{ e_{i}\} $ denote the system of primitive
orthogonal idempotents associated to the above product decomposition of $S$.
As $F$ is a $\Z_{p}$-algebra endomorphism, we know that $ \{
F ( e_{i} )  \}  $ is a system of orthogonal idempotents with 
\begin{equation*}
1=\sum_{i=1}^{n}F ( e_{i} )
\end{equation*}%
and so this system corresponds to a decomposition of the commutative algebra 
$S$ into $n$ components. Since the decomposition of Noetherian commutative
algebras into indecomposable algebras is unique, we must have $F (
e_{i} ) =e_{\pi  ( i ) }$ for some permutation $\pi $ of $\{ 1,\ldots ,n \}$. It will suffice to show that the permutation $%
\pi $ is the identity$.$ Suppose for contradiction that for some $i$, we
have $\pi ( i) =j\neq i$. We know   by definition that 
\begin{equation*}
F ( e_{i} ) \equiv e_{i}^{p}=e_{i}\mod pS
\end{equation*}%
and so 
\begin{equation*}
e_{i}\equiv F ( e_{i} ) \cdot e_{i}\equiv e_{j}\cdot e_{i}=0\mod pS.
\end{equation*}
However, by Theorem 6.7 (page 123) of [CR1] we know that, since $pS$ is
contained in the radical of $S$, $e_{i}\mod pS$ must be a primitive
idempotent of $S/pS$,  and so we have the desired
contradiction.\ \ \ $\square $

\subsection{$\mathbf{Q}_{p}$-$p$-elementary groups}\label{5.1}

Suppose $G$ is a $\Q_{p}$-$p$-elementary group, so that $G$ may be
written as a semi-direct product $C\rtimes P$, where $C$ is a cyclic
normal subgroup of order $s$, which is prime to $p$, and where $P$ is a
$p$-group. We decompose the commutative group ring $\Z_{p}[ C]$ according as the divisors $m$ of $s$%
\begin{equation}\label{eq32}
\Z_{p} [ C ] =\prod\limits_{m}\Z_{p} [ m ],
\end{equation}%
where $\Z_{p} [ m ] $ is the semi-local ring 
\begin{equation*}
\mathbf{Z}_{p} [ m ] =\mathbf{Z} [ \zeta _{m} ] \otimes _{\Z}\Z_{p}
\end{equation*}%
and where $\zeta_{m}$ is a primitive $m$th root of unity in $\Q_{p}^{c}$. For brevity we set $R[m] =R\otimes _{\Z_p}\Z_{p}[m]$, and we note that by Lemma \ref{lem39} $R[m]$
 decomposes as a product of integral domains each of which is  
$p$-adically complete and which possesses  a lift of Frobenius given by the
restriction of the tensor product of the lift of Frobenius of $ R$ and the
Frobenius automorphism of $\Z_p[m]$.

For each $m$ the conjugation action of $P$ on $C$ induces a homomorphism $\alpha _{m}:P\rightarrow {\rm Aut} \langle \zeta _{m} \rangle $ and we
let $H_{m}=\ker (\alpha _{m})$ and $A_{m}={\rm Im}(\alpha _{m})$.

Tensoring the decomposition (\ref{eq32}) with $\otimes_{\Z_p[C]}R[G]  $ affords a decomposition of $R$-algebras%
\begin{equation}\label{eq33}
R[G]=\prod\limits_{m}R[m] \circ P\;
\end{equation}%
where $R[m] \circ P$ denotes the natural twisted group ring. We
shall study the determinant group $\mathrm{Det} ( \GL(R[G]))$
  by studying the various subgroups $\mathrm{Det} ( \GL(R[m]\circ P))$. Note that the twisted group ring $R\left[ m\right] \circ P$ contains the standard group ring $R [ m ]
 [ H_{m} ]$. We therefore have the inclusion map $i:R [ m ]
 [ H_{m} ]  \rightarrow R [ m ] \circ P$. We also
have a restriction map defined by choosing a transversal $\left\{
a_{i}\right\} $ of $P/H_{m}$: This induces a restriction homomorphism 
$$
\mathrm{res}:\GL_{n} ( R [ m ] \circ P ) \rightarrow
\GL_{n\left\vert A_{m}\right\vert } ( R [ m ]  [ H_{m} ]
 ).
$$
 By Proposition \ref{prop13} (b) we know that $\mathrm{Det} ( \GL( R%
 [ m ]  [ H_{m} ] )) =\mathrm{Det} ( R%
 [ m ]  [ H_{m} ] ^{\times } )   $ and so we have
defined the composition:
\begin{equation}\label{eq34}
r_m: {\rm Det} (\GL_{n} ( R [ m ] \circ P )) \rightarrow
{\rm Det}(\GL_{n\left\vert A_{m}\right\vert } ( R [ m ]  [ H_{m} ]
 ))\rightarrow \mathrm{Det} ( R%
 [ m ]  [ H_{m} ] ^{\times } ).  
\end{equation}%
Since for $\pi \in P$, $x\in R [ m ] \circ P^{\times }$, we know
that $ \mathrm{Det} ( \pi x\pi ^{-1} ) =\mathrm{Det} ( x )$, we see that 
\begin{equation*}
r_m: {\rm Det} (\GL_{n} ( R [ m ] \circ P )) \rightarrow
\mathrm{Det} ( R%
 [ m ]  [ H_{m} ] ^{\times } )^{A_m}.
\end{equation*}%
Here $A_{m}$ acts via $\alpha _{m}$ on $R [ m ] $ and by
conjugation on $H_{m}$. From (3.8) on page 69 of [T]  we know that $
r_{m}$ is injective. Note for future reference that for $x\in R [ m]  [ H_{m} ] ^{\times }$, $i ( x ) $ is mapped by
restriction to the diagonal matrix $\mathrm{diag} ( x^{a_{i}} )$;
thus we write ${\rm Det}(x)$ for the usual element
of $\mathrm{Det} ( R [ m ]  [ H_{m} ] ^{\times
} ) $ whereas $\mathrm{Det} ( i ( x )  ) $ denotes an
element of $\mathrm{Det} ( (R [ m ] \circ P)^{\times } )$.
These two determinants are related by the identity
\begin{equation*}
r_{m} ( \mathrm{Det}  (i ( x))))
=\prod\limits_{a\in A_{m}}\mathrm{Det} ( x^{a} ) \overset{\rm defn}{=}%
N_{A_{m}} ( \mathrm{Det} ( x ) ) .
\end{equation*}%
Next we describe $\mathrm{Det} ( (R[m] \circ P)^{\times
})$, and more generally $\mathrm{Det} ( \GL ( R[m]
\circ P ) )$, and the maps $i$ and $r_{m}$ in terms of
character functions. In Lemma \ref{lem46} we shall see that every irreducible
character of $G$ may be written in the form $\mathrm{Ind}_{H_{m}}^{G}(\phi
_{m})$ where $\phi _{m}$ is an abelian character of $ H_{m}$ with the
property that the restriction of $\phi _{m}$ to $C$ has order $m$, for some $%
m$. With this notation the elements ${\rm Det} ( i (
x))$ in $\mathrm{Det} ( R[m] \circ P^{\times
} ) $ are character functions on such $\mathrm{Ind}_{H_{m}}^{G}(\phi
_{m})$ with 
\begin{equation*}
\mathrm{Det} ( i ( x )  ) \left( \mathrm{Ind}%
_{H_{m}}^{G}\phi _{m}\right) =r_{m} ( \mathrm{Det} ( i (
x )  )  )  ( \phi _{m} ) =\prod\limits_{a\in A_{m}}%
\mathrm{Det} ( x^{a} )  ( \phi _{m} ) .
\end{equation*}%
It is also instructive to see the above in the context of K-theory. We then
have induction and restriction maps 
\begin{eqnarray*}
i_{\ast } &:&\Kr_{1} ( R[m][H_m]  )
\leftrightarrows \Kr_{1} ( R[m] \circ P ) :r_{m} \\
i_{\ast } &:&\mathrm{Det} ( R [ m ]  [ H_{m} ]
^{\times } ) \leftrightarrows \mathrm{Det} ( (R [ m ] \circ
P)^{\times } ) :r_{m}.
\end{eqnarray*}%
Similarly we have the corresponding maps on the representation rings%
\begin{equation*}
\Kr_{0} ( \mathbf{Q}_{p}^{c}  [ m%
 ]  [ H_{m} ]  ) \leftrightarrows \Kr_{0} ( \mathbf{Q}_{p}^{c} [ m ] \circ P )
\end{equation*}%
and by Mackey theory the induction map $i=\mathrm{Ind}_{H_{m}}^{P}$ maps 
the ring $
\Kr_{0} ( \mathbf{Q}_{p}^{c}  [ m%
 ]  [ H_{m} ]  )$ onto $\Kr_{0} ( \mathbf{Q}_{p}^{c} [ m ] \circ P )$.
(see page 68 in [T])$.$

Let $I_{H_{m}}$, resp. $I_{P}$, denote the augmentation ideal of $R[ m] [ H_{m}] $, resp. the two sided $R[ m] \circ P$
-ideal generated by $I_{H_{m}}$. The main result of this subsection is to
show:

\begin{thm}\label{thm40}
(a) The map $r_{m}$ defined in (\ref{eq34}) gives an isomorphism 
\begin{equation*}
r_m: {\rm Det} (\GL  ( R [ m ] \circ P ))= \mathrm{%
Det} ( (R [ m ] \circ P)^{\times } ) \xrightarrow{\ r_m\ }
\mathrm{Det} ( R%
 [ m ]  [ H_{m} ] ^{\times } )^{A_m}.
\end{equation*}%

(b) The restriction of $r_{m}$ induces an isomorphism 
\begin{equation*}
\mathrm{Det} ( \GL ( R [ m ] \circ P,I_{P} )  ) =%
\mathrm{Det} ( 1+I_{P} ) \xrightarrow{\ r_m\ }
\mathrm{Det} ( 1+I_{H_{m}} )^{A_{m}}.
\end{equation*}
\end{thm}

Proof. From now on we fix $m$ and therefore drop the index $m$ wherever
possible. We first prove (a). We have seen that $r_{m}$ is injective
on $\mathrm{Det} ( \GL ( R [ m ] \circ P )  ) $; we now show that $r_{m}$ maps $\mathrm{Det} (  ( R [ m ]
\circ P ) ^{\times } ) $ onto $\mathrm{Det} ( R [ m ][ H ] ^{\times } ) ^{A}$.

First we put $\widetilde{P}=P/\left[ H,H\right] $ and we note that from
Theorem 8.3 page 614 of [W2], the ring $R [ m ] \circ \widetilde{P}$
is isomorphic to the ring of $\left\vert A\right\vert \times \left\vert
A\right\vert $ matrices over $ ( R [ m ] [ H^{\rm ab} ]
) ^{A}$. Thus we see that $r$ induces an isomorphism 
\begin{equation}\label{eq35}
\mathrm{Det} ( (R [ m ] \circ \widetilde{P})^{\times } )
\cong  ( R [ m ]  [ H^{\rm ab} ]  ) ^{A\times }.
\end{equation}%
From (\ref{eq17})   and Theorem \ref{thm28} (which trivially extends to products of
rings, since formation of determinants commutes with ring products) and
using (\ref{eq35}), we have   a commutative diagram with exact top row: 
\begin{equation}\label{eq36}
\begin{array}{ccccccc}
1\rightarrow & \mathrm{Det} ( 1+\mathcal{A} ( R [ m ]  [
H ]  )  ) ^{A} & \rightarrow & \mathrm{Det} ( R [ m%
 ]  [ H ] ^{\times } ) ^{A} & \rightarrow &  ( R [
m ]  [ H^{\rm ab} ]  ) ^{A\times } &  \\ 
& \uparrow &  & \uparrow r &  & \uparrow \cong &  \\ 
& \mathrm{Det} ( i ( 1+\mathcal{A} ( R [ m ]  [ H%
 ]  )  )  ) & \subset & \mathrm{Det} ( (R [ m%
 ] \circ P)^{\times } ) & \overset{q}{\twoheadrightarrow } & 
\mathrm{Det} ((R [ m ] \circ \widetilde{P})^{\times } ). & 
\end{array}%
\end{equation}%
It will therefore suffice to show 
\begin{equation*}
r ( \mathrm{Det}(i ( 1+\mathcal{A} ( R[m][H] )))\supseteq \mathrm{Det} ( 1+\mathcal{A}%
 ( R [ m ]  [ H ]  )  ) ^{A}
\end{equation*}%
and this follows from the commutative diagram%
\begin{equation*}
\begin{array}{ccc}
\mathrm{Det}\left( 1+\mathcal{A} ( R [ m ] [H]
\right)  ) & \overset{\nu }{\cong } & \mathrm{\phi }\left( \mathcal{A}%
 ( R[H]) \right) \otimes _{R}R [ m ] \\ 
\downarrow r &  & \downarrow \text{\textrm{Tr}} \\ 
\mathrm{Det}\left( 1+\mathcal{A} ( R [ m ]  [ H ]
\right)  ) ^{A} & \overset{\nu ^{A}}{\cong } & \left( \mathrm{\phi }%
\left( \mathcal{A} ( R [ H ]  ) \right) \otimes _{R}R [
m ] \right) ^{A}.%
\end{array}%
\end{equation*}%
Recall that $F$ is the tensor product of the lift of Frobenius on $R$ with
the Frobenius automorphism of $\Q_p(\zeta_m)/\Q_p$.  Note also that $A$ acts on $R[m]=R\otimes 
\mathbf{Z}_{p} [ m ] $ via the second factor; so, because $G$ is $\Q_{p}$-$p$-elementary, the action of $A$ on $ \langle \chi
 ( G )  \rangle $ factors through $\mathrm{Gal} (\Q_p(\zeta_m)/\Q_p)$;  this guarantees the
actions of $F$ and $A$ commute; hence $\nu $ is an isomorphism of $A$-modules, and this then gives the bottom row in the above diagram.

Since $R[m]$ is a free $R[A]$-module, it follows
that $\phi\left( \mathcal{A} ( R [ H ] )
\right) \otimes _{R}R [ m ] $ is a projective $R [ A ] $%
-module (with diagonal action); and so $\mathrm{Tr}$, and therefore $r$, is
surjective. This then competes the proof of part (a).

The proof of (b) is entirely similar. Firstly, we note that the above
mentioned result of Wall is functorial, so that we have a diagram induced by
the augmentation map $H^{\rm ab}\rightarrow \{ 1 \} $: 
\begin{equation*}
\begin{array}{ccccccccc}
0 & \rightarrow & I_{\widetilde{P}} & \rightarrow & R[m] \circ 
\widetilde{P} & \rightarrow & R[m] \circ A & \rightarrow & 0 \\ 
&  & \downarrow &  & \downarrow &  & \downarrow &  &  \\ 
0 & \rightarrow & {\rm M}_{\left\vert A\right\vert } ( I_{H^{ab}}^{A} ) & 
\rightarrow & {\rm M}_{\left\vert A\right\vert } ( R[m][H^{\rm ab}] ^{A} ) & \rightarrow & {\rm M}_{\left\vert A\right\vert } (
R [ m ] ^{A} ) & \rightarrow & 0%
\end{array}%
\end{equation*}%
Here the central and the right-hand, and hence also the left-hand,
vertical maps are isomorphisms. Secondly, we may restrict (\ref{eq36}) to get the
following commutative diagram, where we use the above to deduce the extreme
right-hand vertical isomorphism:%
\begin{equation*}
\begin{array}{ccccccc}
1\rightarrow & \mathrm{Det} ( 1+\mathcal{A} ( R [ m ] [
H ]  )  ) ^{A} & \rightarrow & \mathrm{Det} (
1+I_{H} ) ^{A} & \rightarrow &  ( 1+I_{H_{m}^{ab}} ) ^{A\times
} &  \\ 
& \uparrow &  & \uparrow r &  & \uparrow \cong &  \\ 
& \mathrm{Det}\left( i ( 1+\mathcal{A} ( R [ m ]  [ H%
 ] )  ) \right) & \rightarrow & \mathrm{Det}\left(
1+I_{P}\right) & \overset{q}{\rightarrow } & \mathrm{Det}\left( 1+I_{%
\widetilde{P}}\right) & \rightarrow 1.\ \ \square%
\end{array}%
\end{equation*}

\begin{prop}\label{prop41}
We have the equality
\begin{equation*}
N_{A} ( \mathrm{Det}\left( 1+I_{H}\right)  ) =\mathrm{Det}\left(
1+I_{H}\right) ^{A}
\end{equation*}%
and so 
\begin{equation*}
\mathrm{Det}\left( 1+I_{P}\right) =  \mathrm{Det}\left( i\left(
1+I_{H}\right) \right)   .
\end{equation*}
\end{prop}

 Proof. For brevity we put $I=I(R[m][H])$, $\overline I=I(R[m][H^{\rm ab}])$.
 Consider the commutative diagram%
\begin{equation*}
\begin{array}{ccccccccc}
0 & \rightarrow & \phi  ( \mathcal{A} ) \otimes R[m] & 
\rightarrow & \mathrm{Det} ( 1+I ) & \rightarrow & 1+\overline{I} & 
\rightarrow & 1 \\ 
&  & {\rm Tr}_{A}\downarrow \uparrow &  & N_{A}\downarrow \uparrow &  & 
N_{A}\downarrow \uparrow &  &  \\ 
0 & \rightarrow & \left( \phi  ( \mathcal{A} ) \otimes R[m] \right) ^{A} & \rightarrow & \mathrm{Det} ( 1+I ) ^{A} & 
\rightarrow & 1+\overline{I}^{A} & \rightarrow & 1.%
\end{array}%
\end{equation*}%
The top row is exact by (\ref{eq16}). The second row is exact on taking $A $-fixed points and noting, as previously that $\phi  ( \mathcal{A}) $   is a cohomologically trivial $A$-module. We have seen previously
that $\mathrm{Tr}_{A}$ is surjective. The surjectivity of the norm map $%
N_{A}:1+\overline{I}\rightarrow 1+\overline{I}^{A}$ is standard and follows
by forming an $A$-stable filtration of $1+I ( R [ H^{\rm ab}%
 ]  )$ whose subquotients are isomorphic to $R/pR$; thus,
when the filtration is extended to $1+\overline{I}=1+I(R[m][H^{\rm ab}])$, then on such subquotients the norm
identifies with the trace $\mod p$ and we argue as before; see page 71
in [T] for details.\ \ \ \ $\square $
\medskip

Finally we show:

\begin{thm}\label{thm42}
Let $S$ be a ring which contains $R$, with both $S$ and $R$ being Noetherian
normal rings which satisfy the hypotheses stated in the Introduction. 
Suppose further that $S$ is a finitely generated free $R$-module and that
the field of fractions of $S$ is Galois over $N$ with Galois group $\Delta $,
which acts on $S$ as $R$-algebra automorphisms, with the property that $S^{\Delta }=R$. Let $ G$ be a finite $\mathbf{Q}_{p}$-$p$-elementary group
and let $\Delta $ act coefficientwise on $\mathrm{Det} ( S[G]
^{\times } )$; then%
\begin{equation*}
\mathrm{Det} ( \GL ( R[G])) =\mathrm{Det}%
 ( R [ G ] ^{\times } )
\end{equation*}%
and 
\begin{equation*}
\mathrm{Det} ( \GL ( S [ G])) ^{\Delta }=%
\mathrm{Det} ( \GL (R[G])) .
\end{equation*}
\end{thm}

Proof. The first equality follows from the decomposition (\ref{eq33}) together with
Theorem \ref{thm40}. We now prove the second equality. By (\ref{eq33}) and 
Theorem \ref{thm40} 
\begin{eqnarray*}
\mathrm{Det} ( \GL ( S [ G])) ^{\Delta }
&=&\oplus _{m}\mathrm{Det} ( S [ m ] \circ P^{\times } )
^{\Delta } \\
&=&\oplus _{m} ( \mathrm{Det} ( S [ m] [H_{m}]^{\times
} ) ^{A_{m}} ) ^{\Delta }.
\end{eqnarray*}%
Recall that $\Delta $ acts via the first term in $S [ m ] \circ
P=S\otimes _{R}(R [ m ] \circ P)$ and that $A_{m}$ acts via the
second term; hence the actions of $\Delta $ and $A_{m}$ commute on $S[m][H_m] =(S\otimes _{R}R[m])[H_m]$; hence we see that 
\begin{equation*}
\mathrm{Det} ( \GL ( S[G]))^{\Delta
}=\oplus _{m}\left(\mathrm{Det} ( S [ m ]  [ H_{m} ]
^{\times })\right) ^{A_{m}\times \Delta }=\oplus _{m}\left( \mathrm{Det}%
 ( S [ m ]  [ H_{m} ] ^{\times } ) ^{\Delta
}\right) ^{A_{m}}
\end{equation*}%
and so by Theorems \ref{thm32} and \ref{thm40} together with (\ref{eq34}) we have equalities 
\begin{eqnarray*}
\mathrm{Det} ( \GL ( S[G]))^{\Delta}
&=&\oplus _{m}\mathrm{Det} ( R[m][H_m])
^{\times } ) ^{A_{m}} \\
&=&\oplus _{m}\mathrm{Det} ( R[m] \circ P^{\times } ) \\
&=&\mathrm{Det} ( \GL(R[G]))
.\;\;\;\;\square
\end{eqnarray*}%
\medskip

\textbf{Application. }\ We conclude this subsection by considering the
implications of the above result for an \textit{arbitrary} finite group $
\widetilde{G}$. We again suppose $S$ to be as Theorem \ref{thm42}. From 12.6 in [S]
we know that we can find an integer $l$ which is prime to $p$, $\mathbf{Q}_{p}$-$p$-elementary subgroups  $H_{i}$ of $\widetilde{G}$, integers $n_{i}$, and $\theta _{i}\in \Kr_{0} ( \Q_{p} [ H_{i} ] )$,
such that 
\begin{equation*}
l\cdot 1_{G}=\sum_{i}n_{i}\cdot {\rm Ind}_{H_{i}}^{\widetilde{G}} (
\theta _{i} ) .
\end{equation*}%
Thus, given $\mathrm{Det} ( x ) \in  \mathrm{Det} ( \GL ( S%
[\widetilde{G}])) ^{\Delta }$, then by the
Frobenius   structure of the module $\mathrm{Det} ( \GL(S [ 
\widetilde{G}])) $ over $\Kr_{0} ( \mathbf{Q}_{p} [
H_{i}]) $ (see \S \ref{4.1}) we have
\begin{equation*}
\mathrm{Det}(x)^{l}=\prod_{i}{\rm Ind}_{H_{i}}^{%
\widetilde{G}} ( \theta _{i}\cdot {\rm Res}_{\widetilde{G}}^{H_{i}} ( \mathrm{Det} ( x))))^{n_{i}}.
\end{equation*}%
However, Theorem \ref{thm42} above implies that

\begin{equation*}
\theta _{i}\cdot {\rm Res}_{\widetilde{G}}^{H_{i}} ( \mathrm{Det}%
(x))\in \mathrm{Det} ( S [ H_{i} ] ^{\times
} ) ^{\Delta }=\mathrm{Det} ( R [ H_{i} ] ^{\times
} ) .
\end{equation*}%
Thus we have shown

\begin{thm}\label{thm43}
For any finite group $\widetilde{G}$ each element in the quotient group 
\begin{equation*}
\mathrm{Det} ( \GL ( S [ \widetilde{G} ]  )  )
^{\Delta }/\mathrm{Det} ( \GL ( R [ \widetilde{G} ]  ))
\end{equation*}%
has finite order which is prime to $p$.
\end{thm}

\subsection{$\mathbf{Q}_{p}$-$\ell $-elementary groups}\label{5.2}

We consider a prime $\ell \neq p$ and a $\mathbf{Q}_{p}$-$\ell $-elementary
group $G$: thus $G$ may be written as $ ( C\times C^{\prime } )
\rtimes L$ where $ C$ is a cyclic $p$-group, $C^{\prime }$ is a
cyclic group of order prime to $p\ell$ and $ L$ is an $\ell $-group. In
this section we show:

\begin{thm}\label{thm44}
If $ G$ is a $\mathbf{Q}_{p}$-$\ell $-elementary group, then $\mathrm{Det}%
( S[G]^{\times } ) ^{\Delta }=\mathrm{Det} ( R%
[G]^{\times } )$ and $\mathrm{Det} ( \GL ( R%
[G]))) =\mathrm{Det} ( R [ G ]^{\times } )$.
\end{thm}

Then, reasoning as in the Application in the previous section, we can
immediately deduce

\begin{thm}\label{thm45}
For an arbitrary finite group $ \widetilde{G}$ each element of the quotient
group $$\mathrm{Det} ( \GL ( S [ \widetilde{G} ]  )  )
^{\Delta }/\mathrm{Det} ( \GL ( R [ \widetilde{G} ]  ))$$
 has finite order which is prime to $\ell$.
\end{thm}

This, together with Theorem \ref{thm43} above, will then
establish Theorem \ref{thm6}.
\medskip

Prior to proving Theorem \ref{thm44}, we first need to recall three preparatory
results:

\begin{lemma}\label{lem46}
Each irreducible character $\chi $ of $G$ can be written in the form ${\rm Ind}_{H}^{G}\phi$, where $\phi$ is an abelian character of a subgroup $H$
which contains $C\times C^{\prime }$.
\end{lemma}

Proof. See 8.2 in [S].\ \ \ \ $\square $

\begin{prop}\label{prop47}
Let $\O$ denote the ring of integers of the finite extension of $\Q_p$ generated by the values of all characters of $G$,
let $\mathfrak{m}$ denote the maximal ideal of $\O$, and let $\P$ denote the $R\otimes _{\mathbf{Z}_{p}}\mathcal{O}$-ideal generated by $%
\mathfrak{m}$. With the notation of the previous lemma, we write $\phi =\phi
^{\prime }\phi _{p}$ where $\phi ^{\prime }$, resp. $\phi _{p}$, has order
prime to $p$, resp. $p$-power order; and we put $\chi ^{\prime }={\rm Ind}_{H}^{G}\phi ^{\prime }$. Then for $r\in \GL(R[G]) $
we have the congruence 
\begin{equation*}
\mathrm{Det} ( r )  ( \chi -\chi ^{\prime } ) \equiv 1%
\mod \P.
\end{equation*}
\end{prop}

Proof. See Theorem 15A on page 66 of [F].\ \ \ \ $\square $

\begin{prop}\label{prop48}
 Put $G^{\prime }=G/C$. Then $\Z_{p} [ G^{\prime } ]
$ is a split maximal $\Z_{p}$-order, that is to say it is a product
of matrix rings 
\begin{equation*}
\mathbf{Z}_{p} [ G^{\prime } ] =\prod_i {\rm M}_{n_{i}} ( \O_{i} )
\end{equation*}%
over (local) rings of integers $\O_i$. Thus we have the equalities
\begin{equation*}
\mathrm{Det} ( \GL( R[ G^{\prime }] ) )
=\prod  ( R\otimes \mathcal{O}_{i} ) ^{\times }=\mathrm{Det} (
R [ G^{\prime } ] ^{\times } )
\end{equation*}%
and
\begin{eqnarray*}
\mathrm{Det} ( S [ G^{\prime } ] ^{\times } ) ^{\Delta }
&=&\prod \mathrm{Det} (  ( S\otimes \mathcal{O}_{i} ) ^{\times
} ) ^{\Delta }=\prod  ( S\otimes \mathcal{O}_{i} ) ^{\times
\Delta } \\
&=&\prod  ( R\otimes \mathcal{O}_{i} ) ^{\times }=\prod \mathrm{%
Det} (  ( R\otimes \mathcal{O}_{i} ) ^{\times } ) \\
&=&\mathrm{Det} ( R [ G^{\prime } ] ^{\times } ) .
\end{eqnarray*}%
\end{prop}

Proof. See Theorem 41.1 and Theorem 41.7 of [Re].\ \ \ \ $\square \medskip $

Proof of Theorem \ref{thm44}. We start by proving the first part. Suppose, as
previously, that $G$ is a $\mathbf{Q}_{p}$-$\ell $-elementary group and that
we are given $\mathrm{Det} ( z ) \in \mathrm{Det} ( S[G]^\times)^\Delta$  and let $z^{\prime }$ denote the image
of $z$ in $S [ G^{\prime } ]$. Then by Proposition \ref{prop48} we know that
we can find $x^{\prime }\in R [ G^{\prime } ] ^{\times }$ with $%
\mathrm{Det} ( x^{\prime } ) =\mathrm{Det} ( z^{\prime } )$; moreover, since the kernel of the map $R[G] \rightarrow R[G^\prime]$
 lies in the radical of $R [ G ]$, we
can of course find $ x\in R[G]^{\times }$ with image $x^{\prime }$ in $R [ G^{\prime } ] ^{\times }$. Thus, to conclude,
it will be sufficient to show that $\mathrm{Det} ( zx^{-1} )$ is in $  
\mathrm{Det} ( R [ G ] ^{\times } )$. However, by
construction, $\mathrm{Det} ( zx^{-1} ) $ is trivial on characters
inflated from $G^{\prime }$, and so by Proposition \ref{prop47} we see that 
\begin{equation*}
\mathrm{Det} ( zx^{-1} )  ( \chi  ) =\mathrm{Det} (
zx^{-1} )  ( \chi -\chi ^{\prime } ) \equiv 1\mod \P
\end{equation*}%
and so we conclude that $\mathrm{Det} ( zx^{-1} ) $ is a pro $p$-element of $\mathrm{Det} ( S [ G^{\prime } ] ^{\times } )^{\Delta }$. But by Theorem \ref{thm43} we know that $\mathrm{Det} (
zx^{-1} ) $ has image in the quotient group $\mathrm{Det} ( S[G ] ^{\times } ) ^{\Delta }/\mathrm{Det} ( R [ G ]
^{\times } ) $ of finite order which is prime to $p$. Therefore we may
deduce that $\mathrm{Det} ( zx^{-1} )$ is in $\mathrm{Det} ( R%
 [ G ] ^{\times } )$.

To prove the second part let $I_{C}$ denote the two-sided $R[G]$-ideal generated by the augmentation ideal $I(R[C])$ of $R[C]$. We then have a commutative diagram 
\begin{equation*}
\begin{array}{ccccccccc}
1 & \rightarrow & \Kr_{1}(R[G] ,I_{C} ) & 
\rightarrow & \Kr_{1} ( R [ G ]  ) & \leftrightarrows & 
\Kr_{1} ( R [ G^{\prime } ]  ) & \rightarrow & 1 \\ 
&  & \downarrow &  & \downarrow &  & \downarrow &  &  \\ 
1 & \rightarrow & \mathrm{Det} ( \GL ( R [ G ]  )
,I_{C} ) & \rightarrow & \mathrm{Det} ( GL ( R [ G ])  ) & \leftrightarrows & \mathrm{Det} ( \GL ( R [
G^{\prime } ]  )  ) & \rightarrow & 1%
\end{array}%
\end{equation*}%
the top row is exact by (\ref{eq4}); this is split exact since the map $G\rightarrow
G^{\prime }$ is split; hence the lower row is split exact. The result then
follows since by Proposition \ref{prop48} and Lemma \ref{lem8} (b) we have
\begin{equation*}
\mathrm{Det} ( \GL(R[G^\prime])) =
\mathrm{Det} ( R[G^\prime]^\times),  \text{ and  }%
\mathrm{Det} ( \GL(R[G ], I_C))=\mathrm{Det} ( \GL(1+I_C)).\text{ \ \ \ }\square
\end{equation*}
\medskip

We can now prove Theorem \ref{thm3} (b) of the Introduction. Let $r\in \GL_{n}(R[G])$. By Brauer induction we can write%
\begin{equation*}
{1}_{G}=\sum\limits_{H}n_{H}\cdot \mathrm{Ind}_{H}^{G} ( \theta
_{H} )
\end{equation*}%
where the $ H$ range over the $\Q_{r}$-elementary subgroups of $G$
for different primes $r$, the $n_{H}$ are integers and $\theta_{H}$ are $\Q_{r}$-characters of $H$. From (\ref{eq29}) we have 
\begin{equation*}
\mathrm{Det}(r) =\sum\limits_{H}n_{H}\mathrm{Ind}_{H}^{G} (
\theta _{H}\cdot {\rm Res}_{G}^{H} ( \mathrm{Det}(r))) .
\end{equation*}%
Now by Theorems \ref{thm42} and \ref{thm44} we know that for each $H$, 
$\mathrm{Det} (
\GL ( R[H])) =\mathrm{Det}(R[H]^\times)$, and so $\theta _{H}\cdot 
{\rm Res}_{G}^{H} ( \mathrm{Det} ( r )  )$ lies in  $\mathrm{Det}(R[H]^\times)$. 
The result then follows since $\mathrm{
Ind}_{H}^{G}$ maps $\mathrm{Det}(R[H]^\times)$
into $\mathrm{Det}(R[G]^\times)$ .\ \ \ \ \
$\square $
\bigskip

\section{Appendix: Torsion determinants}\label{6}

Here we sketch the proof of Theorem \ref{thm26}. The proof is essentially the same as
that given by Wall for the case where $R$ is a $p$-adic ring of integers
and proceeds in three steps. In this section we assume only that $R$ is $p$-adically complete torsion free  
integral domain.
\medskip

{\sl Step 1.} $G$ abelian. We argue by induction on the order of $G$. The result
is trivial if $G=\{1\}$ and this starts the induction. We choose an element $c\in G$ 
of order $p$ and put $\overline{G}=G/\langle c\rangle$.
Consider $x\in ( 1+I(R[G])) _{\rm tors}$.
By the induction hypothesis, we know that, after multiplying $x$ by\ an
element of $G,$ we can assume that the image of $x$ in $R[ \overline{G}] $ is 1, i.e.$\ 1-x\in \left( 1-c\right) R\left[ G\right] ,$ say 
\begin{equation}\label{eq37}
1-x=\left( 1-c\right) y, \text{ with }y\in R[ G] .
\end{equation}

Now choose an abelian character $\theta $ of $G$ which is non-trivial on $c$
and suppose for contradiction that we can find such a $\theta $ with $\theta(x) \neq 1$; by (\ref{eq37}) above, we see that 
$\theta(1-c) $ divides $\theta( 1-x) $ and so $\theta(x) $ must be a primitive $p$-th root of unity. Multiplying $x$ by a
suitable power $c^{a}$ of $c$ we can get $\theta(x)=1$;
moreover, the above work remains valid with this new $x$ but now we also
know that $\theta(y) =0$ and so we can write 
\begin{equation}\label{eq38}
y\in \sum_{g}R ( g-\theta  ( g)) .
\end{equation}%
If $x=1$ we are done; otherwise, by the above, we choose a further abelian
character $\eta $ such that $\eta(x) $ is a root of unity of $p$-power order different from $1.$ Then from (\ref{eq37}) we get 
\begin{equation*}
\eta(1-x) =\eta ( 1-c ) \eta(y) .
\end{equation*}%
Clearly we cannot have $\eta(c) =1,$ for otherwise the
right-hand term vanishes; thus $\eta ( 1-c) $ and $\eta (
1-x)$ must be associates; this, however, is absurd since by (\ref{eq38}) $
\eta(y)$ is a non-unit.
\medskip

{\sl Step 2.} Next we consider the special situation where $G$ contains an
abelian subgroup $H$ of index $p$. We again suppose we have $x\in
1+I(R[G]) $ with the property that $\mathrm{Det}(x)$ is $p$-power torsion, and we show $\mathrm{Det}(x) \in \mathrm{Det}( G) $.

Choosing a transversal of $G/H$ we get a homomorphism 
\begin{equation*}
N:R[ G] ^{\times }\rightarrow \GL_{p}( R[ H]) \overset{\det }{\rightarrow }R[ H] ^{\times }.
\end{equation*}
From 1.4 in [T] we have the restriction map $\mathrm{Res}_{G}^{H}:\mathrm{Det
}( \GL( R[ G])) \rightarrow \mathrm{Det}( \GL( R[ H])) $ with the property that 
\begin{equation}\label{eq39}
\mathrm{Res}_{G}^{H}( \mathrm{Det}( x)) =\mathrm{Det}( N(x)) .
\end{equation}

Choose $\gamma \in G-H$ and let $\alpha $ denote the automorphism of $H$
given by conjugation by $\gamma$; let $T$ denote the submodule of $R[G] $ 
\begin{equation*}
T=\left\{ \sum_{i=1}^{p}\alpha ^{i}(x) \mid x\in R[G]\right\}
\end{equation*}
and let $t:G\rightarrow H$ be the transfer homomorphism. As per Lemma 1.7 on
page 47 of [T] we have

\begin{lemma}\label{lem49}
For $x=\sum_{g\in G}x_{g}g\in R[G]^{\times }$, we have
\begin{equation*}
N(x) \equiv \sum_{g\in G}x_{g}^{p}t(g) \mod T. \ \ \ \square
\end{equation*}  
\end{lemma}

We now complete the proof of this second step. By Step 1 we know that, after
multiplying by a suitable element of $G$ we may assume that $x$ has trivial
image in $R[G^{\rm ab}]$. We now show that for such $x$ with $\mathrm{Det}(x)$ torsion, we must have $\mathrm{Det}(x) =1$.
 To show this we must show that for any non-abelian irreducible
character $\chi $ of $G$,  $\mathrm{Det}(x)(\chi)
=1$. Now, by standard theory, all such $\chi $ are of the form \textrm{Ind}$%
_{H}^{G}(\theta)$ for some abelian character $\theta $ of $H$. Since for a
character $\theta $ of $H$ 
\begin{equation*}
\mathrm{Det}(x)( \text{\textrm{Ind}}_{H}^{G}\theta )
=(\text{\textrm{Res}}_{G}^{H}\mathrm{Det}( x))( \theta)
\end{equation*}%
by (\ref{eq39}) it is enough to show $N(x) =1.$ Note that since $x$ has
trivial image in $R[ G^{\rm ab}]$, 
\begin{equation*}
\mathrm{Det}(N(x)) ( 1_{H} ) =\text{%
(\textrm{Res}}_{G}^{H}\mathrm{Det}( x)) ( 1_{H}) =%
\mathrm{Det}( x) ( \text{\textrm{Ind}}_{H}^{G}1_{H})
=1,
\end{equation*}
and so $N(x) \in 1+I( R[ H] )$. Moreover
by Step 1, since $H$ is abelian, we know that $N(x)=h$ for some 
$h\in H$. Again we write $x=\sum_{g\in G}x_{g}g$; since $x$ has trivial
image in $R[ G^{ab}]$, we deduce that for each $f\in G$ 
\begin{equation*}
\sum_{g\in f\left[ G,G\right] }x_{g}=\left\{ 
\begin{array}{c}
1\text{ if\ }f\in \left[ G,G\right]  \\ 
0\text{ if\ }f\notin \left[ G,G\right]
\end{array} \right.
\end{equation*}
so raising to the $p$-th power we get 
\begin{equation*}
\sum_{g\in f[ G,G] }x_{g}^{p}\equiv \left\{ 
\begin{array}{c}
1 {\mod}p\;\text{ if\ }f\in \left[ G,G\right] \text{ } \\ 
0{\mod}p\;\text{ if\ }f\notin \left[ G,G\right] .%
\end{array} \right.
\end{equation*}%
The restriction of $t$ to $ [ G,G ] $ is trivial and therefore 
\begin{eqnarray*}
h &=&N(x)\equiv \sum_{g\in G}x_{g}^{p}t(g)\mod T \\
&\equiv &\sum_{f\in G/\left[ G,G\right] }t(f) \sum_{g\in f%
\left[ G,G\right] }x_{g}^{p}{\mod}T \\
&\equiv &1 {\mod}( T,p)
\end{eqnarray*}%
and so $h-1\in ( T,p) $. However, $( T,p) \cap
R\subset pR$, and so we conclude that $h=1$, as required.
\medskip

{\sl Step 3.} Consider an arbitrary finite $p$-group $G$ and $\mathrm{Det}(
x) \in \mathrm{Det}( 1+I( R[ G]) )_{\rm tors}$. By Step 1 we may assume that, after multiplying by a suitable
group element, $\widetilde{x}$, the image of $x$ in $R[ G^{\rm ab}]$,
is $1$. We must show that $\mathrm{Det}(x)(\chi)
=1 $ for each non-abelian irreducible character $\chi $ of $G$. By standard
character theory we can write $\chi ={\rm Ind}_{H}^{G}(\theta)$ for some
subgroup $H$ of $G$ of index $p$; we then put $\mathrm{Det} ( x^{\prime
} ) ={\rm Res}_{G}^{H}( \mathrm{Det}( x)) $. Since $\mathrm{Det}( x) $ has trivial image in $R[ G^{\rm ab}]$,
 it follows that 
\begin{equation*}
\mathrm{Det}( x^{\prime }) ({1}_{H}) =\text{%
\textrm{Res}}_{G}^{H}( \mathrm{Det}( x)) ( 
{1}_{H}) =\mathrm{Det}( x) ( \text{\textrm{Ind}%
}_{H}^{G}1_{H}) =1.
\end{equation*}%
So by Proposition \ref{prop13} (a), 
\begin{equation*}
\mathrm{Det}( x^{\prime }) \in \mathrm{Det}( \GL_{p}( R[ H] ,I( R[ H]))) =\mathrm{%
Det}( R[ H] ^{\times })
\end{equation*}%
that is to say we have now shown that we can take $x^{\prime }\in 1+I(R[H]) $ (and not just in $\GL_{p}( R[ H])$.) We then consider the commutative diagram with natural maps: 
\begin{equation*}
\begin{array}{ccccc}
\mathrm{Det}( 1+I( R[ G]) ) & \overset{r}{
\rightarrow } & \mathrm{Det}( 1+I( R[ G/[ H,H] 
]) ) & \overset{q}{\rightarrow } & \mathrm{Det}(
1+I( R[ G^{\rm ab}]) ) \\ 
\downarrow \text{\textrm{Res}} &  & \downarrow \text{\textrm{Res}} &  &  \\ 
\mathrm{Det}( 1+I( R[ H]) ) & \overset{%
r^{\prime }}{\rightarrow } & \mathrm{Det}( 1+I( R[ H^{\rm ab}] ) ) . &  & 
\end{array}%
\end{equation*}%
To conclude: by Step 2 we know that $r( \mathrm{Det}( x)) =1$; hence by the commutative diagram 
\begin{equation*}
r^{\prime }\circ \text{\textrm{Res}}_{G}^{H}( \mathrm{Det}(
x)) =r^{\prime }( \mathrm{Det}( x^{\prime })) =1.
\end{equation*}%
By induction on the group order, we know that $\mathrm{Det}( x^{\prime
}) =1$, and thus 
\begin{equation*}
\mathrm{Det}( x)( \chi) =\mathrm{Det}( x)( \text{\textrm{Ind}}_{H}^{G}\theta) =\mathrm{Res}_{G}^{H}( 
\mathrm{Det}( x))( \theta) =\mathrm{Det}
( x^{\prime })( \theta) =1.\ \ \ \ \ \square
\end{equation*}

\end{document}